\documentclass[11pt,reqno,a4paper]{amsart}

\usepackage{amsthm}
\usepackage{mathtools}
\usepackage{amssymb}
\usepackage{hyperref}
\usepackage{enumitem}
\usepackage{graphicx}
\usepackage{xcolor}
\usepackage{kotex}
\usepackage{comment}

\numberwithin{equation}{section}
\numberwithin{figure}{section}
\allowdisplaybreaks

\theoremstyle{plain}
\newtheorem{thm}{\protect\theoremname}[section]
\newtheorem{lem}[thm]{\protect\lemmaname}
\newtheorem{prop}[thm]{\protect\propositionname}

\theoremstyle{remark}
\newtheorem{rem}{\protect\remarkname}[section]

\providecommand{\theoremname}{Theorem}
\providecommand{\corollaryname}{Corollary}
\providecommand{\lemmaname}{Lemma}
\providecommand{\remarkname}{Remark}
\providecommand{\propositionname}{Proposition}

\begin{document}
\global\long\def\Im{\mathrm{Im}}%
\global\long\def\Re{\mathrm{Re}}%
\global\long\def\s{\sigma}%
\global\long\def\d{\partial}%
\global\long\def\mc#1{\mathcal{#1}}%
\global\long\def\Right{\Rightarrow}%
\global\long\def\Left{\Leftarrow}%
\global\long\def\les{\lesssim}%
\global\long\def\hook{\hookrightarrow}%
\global\long\def\D{\mathbf{D}}%
\global\long\def\rad{\mathrm{rad}}%
\global\long\def\jp#1{\langle#1\rangle}%
\global\long\def\norm#1{\|#1\|}%
\global\long\def\ol#1{\overline{#1}}%
\global\long\def\wt#1{\widehat{#1}}%
\global\long\def\br#1{(#1)}%
\global\long\def\Bb#1{\Big(#1\Big)}%
\global\long\def\bb#1{\big(#1\big)}%
\global\long\def\lr#1{\left(#1\right)}%
\global\long\def\la{\lambda}%
\global\long\def\al{\alpha}%
\global\long\def\be{\beta}%
\global\long\def\ga{\gamma}%
\global\long\def\La{\Lambda}%
\global\long\def\De{\Delta}%
\global\long\def\na{\nabla}%
\global\long\def\fl{\flat}%
\global\long\def\sh{\sharp}%
\global\long\def\calN{\mathcal{N}}%
\global\long\def\avg{\mathrm{avg}}%
\global\long\def\bbR{\mathbf{\mathbb{R}}}%
\global\long\def\bbC{\mathbf{\mathbb{C}}}%
\global\long\def\bbZ{\mathbf{\mathbb{Z}}}%
\global\long\def\bbN{\mathbf{\mathbb{N}}}%
\global\long\def\bbT{\mathbb{T}}%
\global\long\def\bbM{\mathbb{M}}%
\global\long\def\bbP{\mathbb{P}}%
\global\long\def\bbS{\mathbb{S}}%
\global\long\def\bfD{\mathbf{D}}%
\global\long\def\bfL{\mathbf{L}}%
\global\long\def\calF{\mathcal{F}}%
\global\long\def\calH{\mathcal{H}}%
\global\long\def\calL{\mathcal{L}}%
\global\long\def\calE{\mathcal{E}}%
\global\long\def\calO{\mathcal{O}}%
\global\long\def\calR{\mathcal{R}}%
\global\long\def\calD{\mathcal{D}}%
\global\long\def\calI{\mathcal{I}}%
\global\long\def\calT{\mathcal{T}}%
\global\long\def\calY{\mathcal{Y}}%
\global\long\def\calA{\mathcal{A}}%
\global\long\def\calB{\mathcal{B}}%
\global\long\def\calC{\mathcal{C}}%
\global\long\def\calG{\mathcal{G}}%
\global\long\def\calM{\mathcal{M}}%
\global\long\def\calP{\mathcal{P}}%
\global\long\def\calS{\mathcal{S}}%
\global\long\def\calU{\mathcal{U}}%
\global\long\def\Lmb{\Lambda}%
\global\long\def\eps{\varepsilon}%
\global\long\def\lmb{\lambda}%
\global\long\def\gmm{\gamma}%
\global\long\def\rd{\partial}%
\global\long\def\chf{\mathbf{1}}%
\global\long\def\td#1{\widetilde{#1}}%
\global\long\def\sgn{\mathrm{sgn}}%
\global\long\def\blambda{\boldsymbol\lambda}%
\global\long\def\biota{\boldsymbol\iota}%
\global\long\def\Red#1{\textcolor{red}{#1}}%
\global\long\def\NL{\mathrm{NL}}%

\title[Soliton resolution for the Skyrme model]{Stable soliton resolution for the corotational Skyrme model}


\author{Uihyeon Jeong}
\email{uihyeon.jeong@univie.ac.at}
\address{Fakultät für Mathematik, Oskar-Morgenstern-Platz 1, 1090 Wien, Austria}

\begin{abstract}
We study the classical Skyrme model. Under corotational symmetry, we prove that
(i) stationary solutions, or Skyrmions, are unique in every
topological degree; (ii) the Cauchy problem is globally well-posed
in a nonlinear energy space; and (iii) stable soliton resolution holds:
every finite-energy solution scatters to the Skyrmion of the same
degree. Our key observation is that subcritical equivariant Yang--Mills dynamics
emerge near the origin and suppress concentration.
\end{abstract}
\maketitle
\section{Introduction}

\subsection{The classical Skyrme model}
We study the classical Skyrme model, a relativistic field theory of pions
in three space dimensions. Skyrme introduced this model to describe
baryons, including protons and neutrons as spatially localized
configurations of the pion field \cite{Skyrme61}.

A possible route toward such configurations is to introduce a curved target space. In the nonlinear sigma model, pion interactions are described through a field
$\Psi:\bbR^{1+3}\to\mathbb S^3\simeq SU(2)$,
where $\Psi=(\sigma,\pi_1,\pi_2,\pi_3)$ with a scalar $\sigma$ and the three pion components $\pi_1,\pi_2,\pi_3$.
Its dynamics is governed by the wave map equation, the
Euler--Lagrange equation associated with the action
\[
 \mathcal S_{\mathrm{WM}}(\Psi)=\frac12\int_{\bbR^{1+3}}
       \langle\partial_\alpha\Psi,\partial^\alpha\Psi\rangle dtdx,
\]
where $\bbR^{1+3}$ carries the Minkowski metric
$\operatorname{diag}(-1,1,1,1)$ and $\mathbb S^3$ carries the standard
round metric. The corresponding equation conserves the energy
\[
 E_2(\Psi,\Psi_t)=\frac12\int_{\bbR^3}
                    \bigl(|\Psi_t|^2+|\nabla\Psi|^2\bigr)dx.
\]
To describe localized configurations, we impose a fixed vacuum value
at spatial infinity. The spatial domain then admits a one-point
compactification, so that, at each time, the field defines a map
from $\mathbb S^3$ to $\mathbb S^3$. Its topological degree is an
integer preserved by continuous evolution. Skyrme identified this
conserved charge with baryon number \cite{Skyrme61}.

However, the wave map problem is \emph{energy-supercritical}: under
its scaling symmetry
$\Psi_\lambda(t,x)=\Psi(t/\lambda,x/\lambda)$ for $\lambda>0$,
its conserved energy satisfies
\[
 E_2(\Psi_\lambda(t),\partial_t\Psi_\lambda(t))
 =\lambda E_2(\Psi(t),\partial_t\Psi(t)).
\]
This scaling rules out nonconstant finite-energy stationary
solutions and suggests that solutions may concentrate despite
energy conservation. This concentration scenario was confirmed by later finite-time singularity formation results \cite{Shatah88,TurokSpergel1990PhysRevLett}.

Skyrme resolved this obstruction by adding a quartic correction
to the wave map action \cite{Skyrme61}:
\begin{equation}\label{eq:general-action}
 \mathcal S(\Psi)=\mathcal S_{\mathrm{WM}}(\Psi)
 +\frac14\int_{\bbR^{1+3}}
 \langle\partial_\alpha\Psi,\partial^\alpha\Psi\rangle^2
 -\langle\partial_\alpha\Psi,\partial_\beta\Psi\rangle
  \langle\partial^\alpha\Psi,\partial^\beta\Psi\rangle
 dtdx.
\end{equation}
The conserved energy now has the form $E=E_2+E_4$, where $E_4$ is
the quartic contribution. For static fields, the two terms scale in
opposite directions:
\[
 E(\Psi_\lambda,0)
 =\lambda E_2(\Psi,0)+\lambda^{-1}E_4(\Psi,0).
\]
The competition between these two terms introduces a preferred spatial scale and removes the scaling obstruction to stationary solutions.

A localized static field with nonzero degree, if it exists,
cannot be continuously deformed to the vacuum state and is therefore
a \emph{topological soliton}. Such solitons in this model are now called
\emph{Skyrmions} \cite{PakTze1979Annals.Phys}. Skyrme explored
these solitons through approximate calculations based on a
spherically symmetric ansatz \cite{Skyrme61}. It is known that the flow also preserves this equivariant ansatz, corresponding to \emph{corotational symmetry}.

\subsection{Corotational symmetry}
We consider corotational maps of the form 
\[
 \Psi(t,x)=\left(\cos\psi(t,r),\frac{x}{r}\sin\psi(t,r)\right),
 \qquad r=|x|.
\]
In the $SU(2)$ formulation, these are also known
as \emph{hedgehog} solutions. Under the corotational symmetry assumption,
the Euler--Lagrange equation of \eqref{eq:general-action} then reduces to
\begin{equation}\label{eq:skyrme}
 (r^2+2\sin^2\psi)(\psi_{tt}-\psi_{rr})-2r\psi_r
       +\sin(2\psi)(\psi_t^2-\psi_r^2)+\sin(2\psi)
                  \left(1+\frac{\sin^2\psi}{r^2}\right)=0.
\end{equation}
The conserved energy is $E(\psi,\psi_t)=E_2(\psi,\psi_t)+E_4(\psi,\psi_t)$, where
\begin{align*}
 E_2(\psi,\psi_t)&=\frac12\int_0^\infty
       \left[\psi_t^2+\psi_r^2+2\frac{\sin^2\psi}{r^2}\right]r^2dr,\\
 E_4(\psi,\psi_t)&=\int_0^\infty
       \frac{\sin^2\psi}{r^2}\left[\psi_t^2+\psi_r^2
                      +\frac12\frac{\sin^2\psi}{r^2}\right]r^2dr.
\end{align*}
Finite-energy data split into the sectors
\[
 \mathcal E_{\ell,m}=\{(\psi,\psi_t): E(\psi,\psi_t)<\infty,
       \ \psi(0)=\ell\pi,\ \psi(\infty)=m\pi\},\qquad\text{for }\ell,m\in\mathbb Z.
\]
For $(\psi,\psi_t)\in\mathcal E_{\ell,m}$, the map $\Psi$ has
topological degree
\[
 \deg\Psi
 =\frac2\pi\int_0^\infty\sin^2\psi\,\psi_r dr
 =m-\ell.
\]
By the symmetries $\psi\mapsto\psi+k\pi$, $k\in\mathbb Z$, and
$\psi\mapsto-\psi$, we restrict attention to $\mathcal E_{0,n}$ with
$n\ge0$.

\subsection{Main results}
The rigorous mathematical study of \eqref{eq:skyrme} began with Kapitanski and Ladyzhenskaya, who constructed corotational Skyrmions by variational methods
\cite{KapitanskiLadyzhenskaya83}. McLeod and Troy later provided an alternative construction in every topological degree by ODE
methods, and proved uniqueness of corotational Skyrmions in degrees zero and one \cite{McLeodTroy91}. Our first result extends this uniqueness result to all integer degrees.

\begin{thm}[Uniqueness of stationary solutions]\label{thm:stationary}
For every integer $n\ge0$, there exists a unique stationary solution
$(Q_{0,n},0)\in\mathcal E_{0,n}$ to \eqref{eq:skyrme}.
\end{thm}

The first step in passing from the static to the dynamical
problem is to study the Cauchy problem for \eqref{eq:skyrme}. Earlier well-posedness results were established in Sobolev-type
spaces on $\bbR^5$ at the level of
$H^s_{\rm rad}\times H^{s-1}_{\rm rad}$
for the lifted variable $\bigl(\frac{\psi-\phi}{r},\frac{\psi_t}{r}\bigr)$,
where $\phi$ is the fixed reference angle
in \eqref{eq:reference-angle}.

Geba, Nakanishi, and Rajeev proved small-data global well-posedness
in critical Besov-type spaces at $s=5/2$ for data in $\mathcal E_{0,0}$
\cite{GNR12}. For large data, Li proved global regularity in
$H^s_{\rm rad}\times H^{s-1}_{\rm rad}$
for integer $s\ge4$ \cite{Li21}, and Geba and Grillakis
extended this result to real $s>7/2$ \cite{GebaGrillakis18}. 

Building on these works, we develop a Cauchy theory
for all finite-energy data by equipping $\mathcal E_{0,n}$
with the complete metric
\begin{equation}\label{eq:energy-distance}
\begin{aligned}
 d\bigl((\psi,\dot\psi),(\varphi,\dot\varphi)\bigr)
 & \!:=\!\|(\psi-\varphi,\dot\psi-\dot\varphi)\|_{\dot H^1\times L^2}
   +\left\|\frac{\sin\psi\,\psi_r-\sin\varphi\,\varphi_r}{r}\right\|_{L^2}\\
 &\quad+\left\|\frac{\sin\psi\,\dot\psi-\sin\varphi\,\dot\varphi}{r}\right\|_{L^2}
   +\left\|\frac{\sin^2\psi-\sin^2\varphi}{r^2}\right\|_{L^2},
\end{aligned}
\end{equation}
where all norms are taken over $\mathbb R^3$. The metric \eqref{eq:energy-distance} incorporates contributions
from both the quadratic energy $E_2$ and the quartic energy $E_4$,
and induces a topology stronger than the usual
$\dot H^1(\bbR^3)\times L^2(\bbR^3)$ topology. Our next result proves global well-posedness
in the nonlinear energy space $\mathcal E_{0,n}$.

\begin{thm}[Global well-posedness]\label{thm:gwp}
Let $n\ge 0$ be an integer. Then equation \eqref{eq:skyrme} is globally well-posed in $\mathcal E_{0,n}$. More precisely:
\begin{enumerate}[label=(\roman*),widest=iii,leftmargin=*]
\item For any
$(\psi_0,\dot\psi_0)\in\mathcal E_{0,n}$, there exists a unique global
solution $(\psi,\psi_t)\in C(\bbR;\mathcal E_{0,n})$
to \eqref{eq:skyrme} with initial data $(\psi_0,\dot\psi_0)$. 
\item For every $T>0$, the solution map
$(\psi_0,\dot\psi_0)\mapsto(\psi,\psi_t)$ is continuous from
$\mathcal E_{0,n}$ to $C([-T,T];\mathcal E_{0,n})$.
\item The energy is conserved: $E(\psi(t),\psi_t(t))=E(\psi_0,\dot\psi_0)$ for all $t\in\bbR$.
\end{enumerate}
\end{thm}

\begin{rem}[Well-posedness in Sobolev spaces]\label{rem:persistence}
Theorem~\ref{thm:gwp} and Proposition~\ref{prop:persistence}
imply global well-posedness for the lifted variable
$\bigl(\frac{\psi-\phi}{r},\frac{\psi_t}{r}\bigr)$ in
$H^s_{\rm rad}(\bbR^5)\times H^{s-1}_{\rm rad}(\bbR^5)$
whenever $s>5/2$.
This covers all Sobolev regularities above the formal
critical index $s_c=5/2$.
\end{rem}

We next consider the long-time behavior of global solutions.
As noted above, finite-energy data split into disjoint topological
sectors, each containing a Skyrmion. The global Cauchy theory therefore allows us to interpret the long-time behavior of solutions in terms of the dynamics of arbitrary perturbations of these stationary states.

From this perspective, the small-data scattering result of Geba,
Nakanishi, and Rajeev can also be viewed as a nonlinear stability
result for the degree-zero Skyrmion, the vacuum state
\cite{GNR12}; see also \cite{AlejoMaulen22}.
For the nontrivial degree-one Skyrmion, whose uniqueness was
already known, Creek, Donninger, Schlag, and Snelson proved
linear stability \cite{CDSS17}.

These mathematical studies are motivated by the dynamical picture
suggested by the earlier numerical work of Bizo\'n, Chmaj, and
Rostworowski \cite{BCR07}. They observed relaxation to the degree-one Skyrmion under corotational perturbations, with excess energy radiated to spatial infinity. They further predicted
that this relaxation occurs for general smooth finite-energy
corotational data, with the solution scattering to a Skyrmion
of the same degree, i.e., \emph{stable soliton resolution}. Our final result provides a rigorous proof of this prediction.

\begin{thm}[Stable soliton resolution]\label{thm:main}
For every integer $n\ge0$, every solution $(\psi,\psi_t)$ to
\eqref{eq:skyrme} with initial data in $\mathcal E_{0,n}$ scatters
to the Skyrmion $(Q_{0,n},0)$.
More precisely, there is a free wave
$(\psi_L,\partial_t\psi_L)$ satisfying
\begin{equation}\label{eq:linear-angular}
 (\psi_L)_{tt}-(\psi_L)_{rr}-\frac2r(\psi_L)_r
                         +\frac2{r^2}\psi_L=0
\end{equation}
such that, as $t\to+\infty$,
\begin{equation}\label{eq:main decomposition}
(\psi(t),\psi_t(t))=(Q_{0,n},0)
       +(\psi_L(t),\partial_t\psi_L(t))
       +o_{\dot H^1(\bbR^3)\times L^2(\bbR^3)}(1).
\end{equation}
Moreover, all energy in excess of the Skyrmion escapes to spatial infinity:
\begin{equation}\label{eq:excess-energy}
 \lim_{R\to\infty}\lim_{t\to+\infty}
 E(\psi(t),\psi_t(t);r>R)
 =E(\psi(0),\psi_t(0))-E(Q_{0,n},0).
\end{equation}
The analogous backward statement follows by time-reversal symmetry.
\end{thm}

\begin{rem}[Excess energy and convergence]\label{rem:energy-scattering}
One can check that the excess energy in \eqref{eq:excess-energy} equals the conserved energy of the free wave $(\psi_L,\partial_t\psi_L)$:
\[
 E(\psi(0),\psi_t(0))-E(Q_{0,n},0)
 =\frac12\int_0^\infty
 \left[
   (\partial_t\psi_L)^2+(\partial_r\psi_L)^2
   +\frac{2\psi_L^2}{r^2}
 \right]r^2\,dr.
\]
The free wave evolves in
$\dot H^1(\bbR^3)\times L^2(\bbR^3)$.
Its quadratic energy $E_2$ is therefore finite, whereas its
quartic energy $E_4$ need not be finite near the origin. 

In fact, one can strengthen \eqref{eq:main decomposition}
to convergence in the nonlinear energy space $\mathcal E_{0,n}$ by cutting off
the free wave near the origin: as $t\to \infty$,
\[
 (\psi(t),\psi_t(t))
 -(1-\chi)(\psi_L(t),\partial_t\psi_L(t))
 \longrightarrow (Q_{0,n},0)
 \quad\text{in }\mathcal E_{0,n}.
\]
\end{rem}
\subsection{Discussion of the results}
\begin{rem}[Method and novelty]
The main difficulty arises from the supercritical wave map dynamics with respect to the quadratic energy. Scaling suggests that the quartic energy suppresses concentration near the origin. However, this heuristic alone does not yield the a priori estimates needed for a rigorous argument.

Our approach is motivated by the observation in the physics literature \cite{Lukierski80,BizonChmajTabo1999PhysRevD} that the quartic correction
in the Skyrme action is itself a Yang--Mills action and hence, under the equivariant ansatz,
subcritical Yang--Mills dynamics emerge when the quartic energy
dominates (see Section~\ref{subsec:skyrme-yang-mills}).
From this perspective, we observe that Li's nonlinear transformation
\cite{Li21}, introduced to remove derivative nonlinearities,
in fact interpolates between supercritical wave map and subcritical Yang--Mills dynamics
in their five-dimensional formulations (see Remark~\ref{rem:WMYM interpolation}).

These observations translate into subcritical and supercritical
growth of the transformed nonlinearity.
Our new input is a coercive energy bound, uniform in time. Under this estimate, the subcritical growth yields global well-posedness in the energy space and prevents loss of compactness through concentration near the origin, while the supercritical growth provides the exterior decay needed for the channel of energy method developed by Duyckaerts, Kenig,
and Merle \cite{DKM13Soliton, DuyckaertsKenigMerle2014IMRN} and extended to five dimensions by Kenig, Lawrie, and Schlag \cite{KLS14}. We thus obtain compactness, without a profile decomposition, and prove soliton resolution by rigidity.

For uniqueness of stationary solutions, we use McLeod--Troy's inverse-profile
coordinates and positivity estimate \cite{McLeodTroy91}.
Our contribution is to prove that the admissible domain for
stationary inverse profiles is convex and that the transformed
energy density is strictly convex there.
These convexity properties extend uniqueness to every degree.
\end{rem}

\begin{rem}[Soliton resolution]
Soliton resolution predicts that generic solutions of nonlinear
dispersive equations asymptotically decouple into solitons and radiation. This conjecture was motivated by numerical observations for KdV
\cite{ZabuskyKruskal65} and was established by inverse scattering
for various integrable equations; see, for example,
\cite{EckhausSchuur83,Schuur86,BorgheseJenkinsMcLaughlin18,
JenkinsLiuPerrySulem18,ChenLiu21}.
More recent integrable methods yield results for half-wave maps
and the Benjamin--Ono equation
\cite{GerardLenzmann25,GassotGerardMiller26}.

Soliton resolution has also been established for critical
Hamiltonian PDEs, which are typically non-integrable.
Starting with the result of Duyckaerts, Kenig, and Merle \cite{DKM13Soliton}, results have been obtained
for radial energy-critical nonlinear wave equations
\cite{DKM23Soliton,CDKM24,JendrejLawrie23NLW} and equivariant wave maps
\cite{DKMM22,JendrejLawrie25}.
Critical nonlinear Schr\"odinger-type examples include the self-dual
Chern--Simons--Schr\"odinger and Calogero--Moser derivative
nonlinear Schr\"odinger equations \cite{KimKwonOh25,KimKwon26}.
See also the references therein for further results.

In contrast to the preceding results, stable soliton resolution
describes the global asymptotic stability of a unique stationary
solution, with particular emphasis on the nontrivial case.
Earlier results for $\mathbb S^3$-valued maps concern modifications
of the supercritical wave map model underlying \eqref{eq:skyrme}.
Kenig, Lawrie, and Schlag established stable soliton
resolution for $1$-equivariant exterior wave maps \cite{KLS14}, and Kenig,
Lawrie, Liu, and Schlag extended it to all
equivariance classes \cite{KLLSExterior15}. Rodriguez
obtained analogous
results on wormholes \cite{Rodriguez19Corotational,Rodriguez18}. These works rule out the concentration scenario by removing the origin.
Lawrie and Rodriguez resolved the same
problem on the full spatial domain for the Adkins--Nappi model,
a semilinear version of \eqref{eq:skyrme}, under a critical-norm
boundedness assumption \cite{LawrieRodriguez19}. For equivariant wave maps from $\mathbb R\times\mathbb H^2$ to
$\mathbb H^2$, Lawrie, Oh, and Shahshahani obtained an analogous result
for a nonperturbative subset of the harmonic maps
\cite{LawrieOhShahshahani17}. Here, energy criticality leaves open the possibility of concentration
at the origin, but the negative curvature of the target manifold
excludes this scenario.
Despite these differences, these works share a contradiction argument in which profile decomposition produces a minimal
counterexample to soliton resolution.
Our proof instead establishes compactness directly after
subtracting a free wave, without a profile decomposition.
\end{rem}

\begin{rem}[Related work]
Without symmetry assumptions, Esteban and M\"uller proved that
finite-energy Sobolev maps have integer degree
\cite{EstebanMuller92}.
Esteban \cite{Esteban86,Esteban04} and Lin and Yang
\cite{LinYang04} established existence of energy minimizers in
degrees $\pm1$ and conditional existence in other degrees.
For the dynamical problem, Crutchfield and Bell identified a
breakdown of hyperbolicity through linearized analysis
\cite{CrutchfieldBell94}. Wong refined this analysis and established local well-posedness for almost stationary data \cite{Wong11}.

For the two-dimensional equivariant model, nonconcentration of
energy \cite{GebaDaSilva13}, small-data scattering
\cite{GebaNakanishiZhang15}, and large-data global regularity
\cite{Creek13,GebaGrillakis19} have been established.
In five dimensions, McNulty constructed self-similar blow-up
solutions for the strong-field limit equation \cite{McNulty20},
and Chen, McNulty, and Sch\"orkhuber proved their stability
\cite{ChenMcNultySchorkhuber26}.
McNulty extended this stable blow-up result to the full equation \cite{McNulty25}.
\end{rem}

\subsection{The Skyrme correction and Yang--Mills equation}
\label{subsec:skyrme-yang-mills}

The quartic term in the Skyrme action admits an interpretation as a
composite Yang--Mills action, already noted for sphere-valued models
by Lukierski \cite{Lukierski80}. We describe this identification for
maps $\Psi:\mathbb R^{1+d}\to\mathbb S^d$, with $d\geq2$.
Throughout this subsection, repeated spacetime indices are summed,
and indices are raised and lowered using the Minkowski metric
$\operatorname{diag}(-1,1,\ldots,1)$.

Let $\nabla^{\mathbb S^d}$ denote the Levi--Civita connection on $\mathbb S^d$, and define
\[
 \nabla^\Psi\coloneqq\Psi^*\nabla^{\mathbb S^d}
\]
on the pullback tangent bundle $\Psi^*T\mathbb S^d$.
Since $\nabla^{\mathbb S^d}$ preserves
the standard metric and orientation, its pullback $\nabla^\Psi$
is also an $SO(d)$-connection.
Under the curvature convention for $\mathbb S^d$:
\[
 R^{\mathbb S^d}(X,Y)=XY^T-YX^T,
 \qquad X,Y\in T_p\mathbb S^d \subset\mathbb R^{d+1},
\]
the pullback connection has curvature
\begin{equation*}
 F_{\alpha\beta}
 =[\nabla^\Psi_\alpha,\nabla^\Psi_\beta]=R^{\mathbb S^d}(\partial_\alpha\Psi,\partial_\beta\Psi)=\partial_\alpha\Psi\,(\partial_\beta\Psi)^T
   -\partial_\beta\Psi\,(\partial_\alpha\Psi)^T.
\end{equation*}
The trace identity $ \mathrm{tr}(ab^Tcd^T)=\langle b,c\rangle\mathrm{tr}(ad^T)
                    =\langle b,c\rangle\langle d,a\rangle$
yields
\begin{align*}
 -\frac12\mathrm{tr}(F_{\alpha\beta}F^{\alpha\beta})
 &=
 \langle\partial_\alpha\Psi,\partial^\alpha\Psi\rangle^2
 -
 \langle\partial_\alpha\Psi,\partial_\beta\Psi\rangle
 \langle\partial^\alpha\Psi,\partial^\beta\Psi\rangle.
\end{align*}
With $\mathcal S_4$ denoting the quartic term in
\eqref{eq:general-action}, we have the following identity
\begin{equation*}
 \mathcal S_4[\Psi]
 =\frac{1}{2}\mathcal S_{\mathrm{YM}}[A(\Psi)],\quad \mathcal S_{\mathrm{YM}}[A]
 \coloneqq -\frac{1}{4}\int_{\mathbb R^{1+d}}
 \operatorname{tr}(F_{\alpha\beta}F^{\alpha\beta})dtdx,
\end{equation*}
where $A(\Psi)$ denotes a gauge potential of $\nabla^\Psi$.

In particular, under our corotational ansatz, parallel transport
of the standard tangent frame at the north pole along meridians yields
\[
 A_0=0,\qquad
 A_k^{ij}
 =\frac{1-\cos\psi(t,r)}{r^2}
 \bigl(x_i\delta_k^{\,j}-x_j\delta_k^{\,i}\bigr),
 \qquad i,j,k=1,\ldots,d.
\]
With $v=\cos\psi$, the connection takes the standard
$SO(d)$-equivariant Yang--Mills form, as observed for $d=3$
by Bizo\'n, Chmaj, and Tabor \cite{BizonChmajTabo1999PhysRevD}.
The corresponding equation is
\begin{equation}\label{eq:strong-field-yang-mills}
 v_{tt}-v_{rr}-\frac{d-3}{r}v_r
 +\frac{d-2}{r^2}v(v^2-1)=0.
\end{equation}

Now, \eqref{eq:skyrme} can be written as follows: for $v=\cos\psi$,
\begin{equation}\label{eq:WMYM 5D variables}
      \psi_{tt}-\psi_{rr}-\frac{2}{r}\psi_r +\frac{\sin(2\psi)}{r^2} = \frac{2\sin\psi}{r^2}\left( v_{tt}-v_{rr}+\frac{v(v^2-1)}{r^2} \right).
\end{equation}
For the five-dimensional formulation, we define
\[g\coloneqq\frac{\psi}{r},\qquad
 w\coloneqq\frac{1-\cos\psi}{r^2},\qquad
 \Box_5\coloneqq\partial_t^2-\partial_r^2-\frac4r\partial_r.\]
Then, we can also rewrite \eqref{eq:skyrme} as
\begin{equation}\label{eq:WMYM reform}
\Box_5g-\frac{2rg-\sin(2rg)}{r^3}
 =-2\frac{\sin \psi}{r}
 \left(\Box_5w-3w^2+r^2w^3\right).
\end{equation}

These observations suggest the following mechanism by which
the quartic energy $E_4$ suppresses concentration at the origin. The main difficulty arises from the supercritical wave map dynamics associated with the quadratic energy.
If the quadratic energy density remains larger than the quartic one
near the origin, $g$ stays uniformly bounded, and so finite-time blow-up
cannot occur. Any possible concentration scenario at the origin would
therefore have to involve unbounded $\sin\psi/r$. In this regime,
the quartic energy dominates, and the subcritical Yang--Mills
dynamics are expected to suppress concentration. We further clarify how the dynamics depend on the size of
$\sin\psi/r$ in Remark~\ref{rem:WMYM interpolation} of Section~3.

\subsection{Notation}
Let $\phi$ be a smooth reference angle such that
\begin{equation}\label{eq:reference-angle}
 \phi(r)=0\quad\text{for }0\le r\le1,\qquad
 \phi(r)=n\pi\quad\text{for }r\ge2.
\end{equation}
For $n>0$, we choose $\phi'(r)>0$ on $1<r<2$; for $n=0$, we take
$\phi\equiv0$. Thus $\phi$ varies only on the fixed annulus
$1<r<2$, and $\boldsymbol\phi=(\phi,0)$.

For nonnegative quantities, $A\lesssim B$ means $A\le CB$ for a
constant $C>0$, which may change from line to line; $A\sim B$
means $A\lesssim B\lesssim A$.
Subscripts indicate allowed dependence; for example,
$A\lesssim_E B$ allows dependence on the energy $E$.
We suppress dependence on the fixed reference angle $\phi$
throughout. Constants in uniform estimates are independent of the
initial time and of the cutoff radii unless indicated otherwise.

We fix a smooth radial cutoff $\chi$ satisfying
$0\le\chi\le1$ and $|\chi'|\lesssim1$, with $\chi(r)=1$
for $r\le1$ and $\chi(r)=0$ for $r\ge2$. We also denote $\chi_R(r)\coloneqq\chi(r/R)$. The symbol $\mathbf1_A$ denotes the characteristic function of a set
$A$.

We use bold notation for vectors in $\bbR^2$ and
$\bbR^2$-valued functions, writing $\boldsymbol f=(f,\dot f)$; along a solution,
$\boldsymbol f(t)=(f(t),\partial_t f(t))$.
For a stationary profile, we write $\boldsymbol Q=(Q,0)$ with $Q=Q_{0,n}$.
Multiplication by a spatial cutoff acts on both components.

The nonlinear transformation in Section~\ref{sec:notation} reduces \eqref{eq:skyrme} to a radial semilinear wave equation on $\bbR^{1+5}$.
We use the following five-dimensional notation for the transformed equation.
We suppress the factor $|\mathbb S^{d-1}|$ in radial integrals and
norms on $\bbR^d$. Unless another dimension is displayed, spatial
norms are on $\bbR^5$; in particular, $\|\cdot\|_p$ means
$\|\cdot\|_{L^p(\bbR^5)}$, and $L^p_I L^q_x$ means
$L^p(I;L^q(\bbR^5))$. Norms restricted to an annulus or a ball use
the same spatial measure. The corresponding energy space for the transformed equation is
\[
 \mathcal H=\dot H^1_{\rm rad}(\bbR^5)\times L^2_{\rm rad}(\bbR^5),
 \qquad
 \|\boldsymbol f\|_{\mathcal H}^2
 =\int_0^\infty(f_r^2+\dot f^{\,2})r^4dr.
\]
A restriction such as $\mathcal H(r>R)$ restricts the defining integral
to $r>R$. We also use the local norm
\begin{equation}\label{eq:local-norm}
 \|\boldsymbol f\|_{\mathcal H_B}^2
 :=\int_0^B(f_r^2+\dot f^{\,2})r^4dr+\int_0^B f^2r^2dr,
 \qquad B>0.
\end{equation}

We set $\Delta_5=\partial_r^2+4r^{-1}\partial_r$. We denote by $S(t)$ the free wave propagator on $\mathcal H$
associated with $W_{tt}-\Delta_5W=0$, so that
$S(t)\boldsymbol W(0)=\boldsymbol W(t)$.

\subsection{Outline of the proof}
Section~\ref{sec:stationary} is devoted to proving uniqueness of stationary solutions. In
Section~\ref{sec:gwp}, we introduce our nonlinear transform and prove global well-posedness. Sections~\ref{sec:radiation} and~\ref{sec:rigidity} are devoted to proving soliton resolution. In Section~\ref{sec:radiation}, we construct the exterior radiation and establish the precompactness of the remainder trajectory obtained by subtracting the radiation component from the solution. Section~\ref{sec:rigidity} is devoted to rigidity of the precompact trajectory, and hence soliton resolution.
We now outline the arguments of the four main sections.

\medskip\noindent
\emph{1. Uniqueness of stationary solutions.}
We adopt the inverse-profile coordinates used in
\cite{McLeodTroy91}. For a nonconstant stationary profile $Q$,
we set
\[
 r=Q^{-1}(\theta),\qquad
 \boldsymbol{x}=(x,\dot{x})=(\log r(\theta),x'(\theta)).
\]
The static energy then takes the form $E(\boldsymbol{Q})=\int_0^{n\pi}L(\theta,\boldsymbol{x})d\theta$, where
\[
 L(\theta,\boldsymbol{x}) \coloneqq e^x\left(\frac1{2\dot{x}}+\dot{x}\sin^2\theta\right)
       +e^{-x}\left(\frac{\sin^2\theta}{\dot{x}}
                          +\frac{\dot{x}\sin^4\theta}{2}\right).
\]
McLeod--Troy's positivity estimate \eqref{eq:static-Q} places every stationary
inverse profile in the convex region
\[
 D_\theta=\left\{\boldsymbol{x}=(x,\dot{x})\in\bbR\times(0,\infty):
 1/\dot{x}>|\sin\theta| /\sqrt2 \right\}.
\]
We prove that $L(\theta,\cdot)$ is strictly convex on $D_\theta$, which yields uniqueness.

\medskip\noindent
\emph{2. Global well-posedness.} We first introduce a simplified version of Li's transformation \cite{Li21}: for $\phi$ in \eqref{eq:reference-angle}, we define $\boldsymbol{u}\coloneqq (u,u_t)$, where
\[u(t,r)=\frac1{r^2}\int_{\phi(r)}^{\psi(t,r)}
                    \sqrt{r^2+2\sin^2s}ds.\]
By Proposition~\ref{prop:nonlinear-transform},
$\boldsymbol{\psi}\mapsto\boldsymbol{u}$ is a homeomorphism
from $\mathcal E_{0,n}$ onto $\mathcal H$, and
\[
 \|\boldsymbol{u}\|_{\mathcal H}
       \lesssim 1+\sqrt{E(\boldsymbol{\psi})}.
\]
Moreover, this transform not only removes the derivative nonlinearity, but also yields the quadratic and cubic bounds near the origin and infinity: 
\[
 u_{tt}-\Delta_5u=F(r,u),\quad |F(r,u)|\lesssim\min\{|u|^2,|u|^3\}
 \quad\text{for }0<r\le1\text{ or }r\ge2,
\]
with the corresponding difference estimates. On $1<r<2$, $F$ is smooth and globally Lipschitz in $u$. Hence, the Strichartz estimate and a contraction argument prove local well-posedness with a lifespan depending only on the initial $\mathcal H$ norm. Once energy conservation is justified by smooth approximation,
the energy bound yields global well-posedness for $\boldsymbol u$
in $\mathcal H$. We thus obtain Theorem~\ref{thm:gwp} through the homeomorphism
$\boldsymbol u\mapsto\boldsymbol\psi$.

\medskip\noindent
\emph{3. Soliton resolution.}
The proof of soliton resolution for $\boldsymbol u$ consists of compactness and rigidity, established in Propositions~\ref{prop:compactness} and~\ref{prop:rigidity}, respectively. 

First, Proposition~\ref{prop:compactness} constructs a free wave
$\boldsymbol u_L$ such that
\[
 K\coloneqq
 \overline{\{\boldsymbol u(t)-\boldsymbol u_L(t):t\ge0\}}^{\mathcal H}
 \quad\text{is compact}.
\]
Hence, every accumulation point of $\boldsymbol u(t)-\boldsymbol u_L(t)$ as $t\to+\infty$ generates a complete solution $\boldsymbol u_*(s)$ such that $\{\boldsymbol{u}_*(s):s\in \mathbb{R}\}\subset K$. Proposition~\ref{prop:rigidity} implies $\boldsymbol u_*(s)\equiv\boldsymbol u_Q$. By compactness of $K$, we conclude that
\[
 \lim_{t\to\infty}
 \|\boldsymbol u(t)-\boldsymbol u_L(t)-\boldsymbol u_Q\|_{\mathcal H}=0.
\]
We then deduce Theorem~\ref{thm:main}, with $\psi_L(t,r)=ru_L(t,r)$ solving \eqref{eq:linear-angular}.

By the uniform energy bound, the radial Sobolev inequality, and
the cubic bound on $F$, we have exterior bounds for $u(r)$ and $F$: for $r>2$ and $R>2$, 
\[
 \sup_{t\in\bbR}|u(t,r)|
 \lesssim_E r^{-3/2},\quad
 \sup_{t\in\bbR}
 \|\mathbf1_{r>R}F(r,u(t))\|_{L^2(\bbR^5)}
 \lesssim_E R^{-2}.
\]
These exterior estimates allow us to apply the standard channel of energy
method
\cite{DKM13Soliton,DuyckaertsKenigMerle2014IMRN,KLS14,
KLLSExterior15,LawrieRodriguez19,DuyckaertsYang21}
in the exterior analysis.
In particular, the proof of Proposition~\ref{prop:rigidity}
follows the argument in
\cite[Sections~5.3--5.5]{LawrieRodriguez19}.
We therefore focus on the proof of
Proposition~\ref{prop:compactness}.

On each fixed ball $\{r<B\}$, the strong Huygens principle for the free wave $\boldsymbol u_L$
reduces the problem to precompactness of
$\{\boldsymbol u(t):t\ge0\}$ in $\mathcal H_B$.
We prove this local compactness in Proposition~\ref{prop:local-compactness}.
For $t>L>B+2$, the Duhamel formula yields
\[
 \boldsymbol u(t)=S(t)\boldsymbol u(0)
   +\int_0^{t-L}\!\!S(t-s)(0,F(r,u(s)))d s
   +\!\int_{t-L}^t\!\!S(t-s)(0,F(r,u(s)))d s.
\]
The first term tends to zero in $\mathcal H_B$ as $t\to+\infty$
by the strong Huygens principle.
For the middle integral term, the strong Huygens principle and the exterior estimate for $F$ yield
\[
 \begin{aligned}
 \left\|\int_0^{t-L}S(t-s)(0,F(r,u(s)))ds\right\|_{\mathcal H_B}
 \lesssim_E(L-B)^{-1}.
 \end{aligned}
\]
Hence, it suffices to prove that the last integral ranges over
a precompact subset of $\mathcal H$ for each fixed $L>B+2$.
Continuity of the Duhamel integral from
$L^1_{[0,L]}L^2_x$ to $\mathcal H$ and the exterior estimate for $F$ reduce this to proving that, for every fixed $R>2$,
\[
 \{\mathbf1_{r<R} F(r,u(\tau+\cdot,r)):\tau\ge0\}
 \quad\text{is precompact in }L^1_{[0,L]}L^2_x.
\]
On $\{r<a\}$ for $0<a<1$, the quadratic bound on $F$ implies
\[
 \sup_{\tau\ge0}
 \|\mathbf1_{r<a}F(r,u)\|_{L^1_{[\tau,\tau+L]}L^2_x}
 \lesssim_{E,L}a^{1/6}.
\]
Precompactness on each fixed annulus $a<r<R$, together with
the bound above, yields local compactness.

Proposition~\ref{prop:radiation} establishes the required
exterior tightness using the channel of energy argument.
Following \cite[Proposition~3.4, Steps~2 and~3]{DuyckaertsYang21},
we first construct a free wave $\boldsymbol u_L$ satisfying
\begin{equation}\label{eq:strategy convergence}
\|\boldsymbol u(t)-\boldsymbol u_L(t)\|_{\mathcal H(r>t+A)}
 \to0
 \qquad\text{as }t\to+\infty
\end{equation}
for every fixed $A\in\bbR$. Then we conclude the proof by showing
\[
 \limsup_{T\to\infty}
 \|\boldsymbol u(T)-\boldsymbol u_L(T)\|_{\mathcal H(r>R)}^2
 \lesssim_E R^{-1},
 \qquad\text{for }R\ge4.
\]
Let $\boldsymbol{w}=\boldsymbol{u}-\boldsymbol{u}_L$ and $T \ge 0$. For the channel argument on $r>R$, we consider
$S(s)\boldsymbol w(T)$. The exterior estimate for $F$ implies
\begin{equation}\label{eq:strategy comparison}
\|S(s)\boldsymbol w(T)-\boldsymbol w(T+s)\|_{\mathcal H(r>R+|s|)}
 \lesssim_E R^{-1}
 \qquad\text{for }s\ge-T.
\end{equation}
\eqref{eq:strategy convergence} with $t=T+s$ and $A=R-T$ yields $\|S(s)\boldsymbol w(T)\|_{\mathcal H(r>R+|s|)} = O_E(R^{-1})$ as $s\to+\infty$.
On the other hand, monotonicity of the exterior free wave energy implies that the same quantity is bounded by $\|S(-T)\boldsymbol w(T)\|_{\mathcal H(r>R+T)}$ as $s\to -\infty$. The channel inequality in Lemma~\ref{lem:channel} and \eqref{eq:strategy comparison} at $s=-T$ then imply
\[
 \|\pi_R^\perp\boldsymbol w(T)\|_{\mathcal H(r>R)}
 \lesssim_E R^{-1}
 +\|\boldsymbol w(0)\|_{\mathcal H(r>R+T)},
\]
where $\pi_R^\perp$ denotes the orthogonal projection onto
$\operatorname{span}\{(r^{-3},0),(0,r^{-3})\}^{\perp}$.
For each fixed $R$, we have $\|\boldsymbol w(0)\|_{\mathcal H(r>R+T)}=o_{E}(1)$ as $T\to \infty$. We iterate \eqref{eq:dyadic-increments} over dyadic radii,
starting from a fixed radius $R_0=4$ where the coefficients of $\pi_{R_0}\boldsymbol w(T)$
are $O_E(1)$. Their growth is then sufficiently slow for
\eqref{eq:projection-split} to yield the desired estimate.
Propositions~\ref{prop:local-compactness} and~\ref{prop:radiation} imply
Proposition~\ref{prop:compactness}. 

\vspace{5pt}
\noindent\mbox{\textbf{Acknowledgements.~}} 
The author would like to thank Sung-Jin Oh for valuable discussions,
as well as Pierre Germain and Roland Donninger for helpful comments
on the manuscript.
\section{Uniqueness of stationary solutions}\label{sec:stationary}
In this section, we prove Theorem~\ref{thm:stationary}.
We express the static energy in inverse-profile coordinates
and establish strict convexity of its density on the convex
domain admissible for stationary profiles. The static energy is
\begin{equation}\label{eq:static-energy-functional}
 E(\boldsymbol{Q})=\int_0^\infty\left[
 \left(\frac{r^2}{2}+\sin^2Q\right)(Q')^2
 +\sin^2Q+\frac{\sin^4Q}{2r^2}\right]d r
\end{equation}
and its Euler--Lagrange equation is
\begin{equation}\label{eq:static-ode}
 (r^2+2\sin^2Q)Q''+2rQ'+\sin(2Q)(Q')^2
 -\sin(2Q)\left(1+\frac{\sin^2Q}{r^2}\right)=0.
\end{equation}
Every finite-energy stationary distributional solution to \eqref{eq:static-ode} is smooth on $(0,\infty)$ by one-dimensional elliptic regularity on compact annuli.
We recall the results of McLeod--Troy \cite[Theorems~2.3--2.5, Lemmas~6.3 and 6.4]{McLeodTroy91} in our normalization, with
their radial coordinate equal to $2r$.

\begin{lem}[Stationary profiles, \cite{McLeodTroy91}]\label{lem:MT-inputs}
For every integer $n\ge0$, there exists a stationary solution
$\boldsymbol Q=(Q,0)\in\mathcal E_{0,n}$. The only stationary profile in
$\mathcal E_{0,0}$ is $Q\equiv0$. For $n\ge1$, every finite-energy
stationary profile in $\mathcal E_{0,n}$ has endpoints
\begin{equation}\label{eq:static-endpoints}
 Q(0)=0,\qquad Q(\infty)=n\pi,\qquad\text{for }n\ge1,
\end{equation}
and satisfies the following properties.
\begin{itemize}
\item (Monotonicity) We have $Q'>0$ and for $r>0$,
\begin{equation}\label{eq:static-Q}
 (r^2+2\sin^2Q)(Q')^2 >2\sin^2Q+\frac{\sin^4Q}{r^2}.
\end{equation}
\item (Asymptotics) There are $\alpha,\beta>0$ such that
\begin{align}
 Q(r)&=\alpha r+O(r^3),& Q'(r)&=\alpha+O(r^2)
 && \text{as }r\downarrow0,\label{eq:static-origin}\\
 n\pi-Q(r)&=\beta r^{-2} +o(r^{-2}),& Q'(r)&=2\beta r^{-3}+o(r^{-3})
 &&\text{as }r\to\infty.\label{eq:static-infty}
\end{align}
\end{itemize}

\end{lem}
\medskip\noindent

For every stationary solution $\boldsymbol Q=(Q,0)\in\mathcal E_{0,n}$
with $n\ge1$, positivity of $Q'$ yields the inverse profile
$r=Q^{-1}(\theta)$ for $\theta\in(0,n\pi)$.
Let $\boldsymbol{x}=(x,\dot{x})$ by
\[
 x(\theta)=\log r(\theta),\qquad
 \dot{x}=x'(\theta)=\frac1{rQ'(r)}>0,\qquad a(\theta)=\sin^2\theta.
\]
In terms of $\boldsymbol{x}$, \eqref{eq:static-energy-functional} becomes
\begin{equation}\label{eq:static-inverseenergy}
 E(\boldsymbol{Q})=\int_0^{n\pi}L(\theta,\boldsymbol{x})d\theta,
 \quad L\coloneqq L_++L_-,
\end{equation}
where $L_+$ and $L_-$ correspond to $E_2(\boldsymbol{Q})$ and
$E_4(\boldsymbol{Q})$, respectively:
\[
 L_+\coloneqq e^x\left(\frac1{2\dot{x}}+a\dot{x}\right),\quad
 L_-\coloneqq e^{-x}\left(\frac{a}{\dot{x}}
                         +\frac{a^2\dot{x}}{2}\right).
\]
Since $r=e^x$ and $rQ'=1/\dot{x}$, \eqref{eq:static-Q} implies
\begin{equation}\label{eq:static-cone}
 \frac1{\dot{x}^2}>a\frac{2+ae^{-2x}}{1+2ae^{-2x}}
 =\frac a2+\frac{3a}{2(1+2ae^{-2x})}\ge\frac a2.
\end{equation}
Thus the inverse profile of every stationary solution
lies in the convex set
\begin{equation}\label{eq:static-domain}
 D_\theta\coloneqq
 \{\boldsymbol{x}:x\in\mathbb R,\ 1/\dot{x}>\sqrt{a/2}\}.
\end{equation}

Now, we prove Theorem~\ref{thm:stationary} by showing strict convexity of $L$ on $D_\theta$.
\begin{lem}\label{lem:static-convex}
For each $\theta\in(0,n\pi)$, $L(\theta,\cdot)$ is strictly
convex on $D_\theta$.
\end{lem}
\begin{proof}
We compute the Hessian determinants with respect to $\boldsymbol{x}$:
\[
 \det D^2L_+
 =e^{2x}\left(\frac1{4\dot{x}^4}
             +\frac{2a}{\dot{x}^2}-a^2\right),\quad
 \det D^2L_-
 =e^{-2x}\left(\frac{a^2}{\dot{x}^4}
             +\frac{2a^3}{\dot{x}^2}-\frac{a^4}{4}\right).
\]
$e^{-2x}\det D^2L_+$ is strictly
increasing in $\dot{x}^{-2}$ on $[0,\infty)$ for every $a\ge0$, with value $a^2/16$
at $\dot{x}^{-2}=a/2$. The inequalities $\dot{x}^{-2}>a/2$ and
$(L_+)_{xx}=L_+>0$ therefore imply that $D^2L_+ >0$.
If $a>0$, $e^{2x}\det D^2L_-$ is strictly increasing as well, with
value $a^4$ at $\dot{x}^{-2}=a/2$, so $D^2L_- >0$. If $a=0$, then $L_-=0$. In either case,
$D^2L=D^2L_++D^2L_- > 0$ on $D_{\theta}$.
\end{proof}

\begin{proof}[Proof of Theorem~\ref{thm:stationary}]
Lemma~\ref{lem:MT-inputs} supplies existence and the classification
in degree zero. For $n\ge1$, let $\boldsymbol Q_i=(Q_i,0)$, $i=1,2$, be stationary solutions in
$\mathcal E_{0,n}$,
with inverse coordinates $\boldsymbol{x}_1,\boldsymbol{x}_2$.
By Lemma~\ref{lem:static-convex},
\begin{equation}\label{eq:convexity another form}
(DL(\theta, \boldsymbol{x}_2)-DL(\theta, \boldsymbol{x}_1)) \cdot (\boldsymbol{x}_2-\boldsymbol{x}_1) \ge 0,
\end{equation}
with equality only if $\boldsymbol{x}_1 = \boldsymbol{x}_2$.
Each $\boldsymbol{x}_i$ satisfies the Euler--Lagrange equation
\begin{equation}\label{eq:static-EL}
 \frac{ d}{ d\theta} \partial_{\dot{x}} L(\theta,\boldsymbol{x}_i)
 =\partial_x L(\theta,\boldsymbol{x}_i).
\end{equation}
Integration by parts using \eqref{eq:convexity another form} and
\eqref{eq:static-EL} yields
\begin{align*}
0&\le \int_0^{n\pi}
 (DL(\theta, \boldsymbol{x}_2)-DL(\theta, \boldsymbol{x}_1))
 \cdot (\boldsymbol{x}_2-\boldsymbol{x}_1)d\theta \\
 &=[(\partial_{\dot{x}}L(\theta,\boldsymbol{x}_2)
          -\partial_{\dot{x}}L(\theta,\boldsymbol{x}_1))(x_2-x_1)
      ]_{0+}^{n\pi-}=0.
\end{align*}
The endpoint terms vanish by \eqref{eq:static-origin},
\eqref{eq:static-infty}, and the following identity
\[
 \partial_{\dot{x}}L(\theta, \boldsymbol{x})=r\left(a-\frac1{2\dot{x}^2}\right)
 +\frac1r\left(\frac{a^2}{2}-\frac{a}{\dot{x}^2}\right).
\]
Strict convexity implies $\boldsymbol{x}_1=\boldsymbol{x}_2$, so
$Q_1=Q_2$. This proves Theorem~\ref{thm:stationary}.
\end{proof}
\section{Global well-posedness}
\label{sec:gwp}
In this section, we prove Theorem~\ref{thm:gwp} by transforming
\eqref{eq:skyrme} into a semilinear wave equation on $\bbR^5$.
We first establish the homeomorphism
$\mathcal E_{0,n}\to\mathcal H=\dot{H}^1_{\rm rad}(\bbR^5)\times L^2_{\rm rad}(\bbR^5)$ and the estimates for the
transformed nonlinearity, then develop the Cauchy theory
in $\mathcal H$.

\subsection{Nonlinear transform}\label{sec:notation}
 For $\boldsymbol\psi\in\mathcal E_{0,n}$ and $r>0$, we define
\begin{equation}\label{eq:transform}
 \boldsymbol{u}
 =\left(\frac{H(r,\psi)-H(r,\phi(r))}{r^2},
         \frac{h(r,\psi)\dot\psi}{r^2}\right),
\end{equation}
where $\phi$ be the reference angle in \eqref{eq:reference-angle} and
\[
 h(r,\psi)=\sqrt{r^2+2\sin^2\psi},\qquad
 H(r,\psi)=\int_0^\psi h(r,s)d s.
\]
For each $r>0$, $H(r,\cdot)$ is strictly increasing onto $\bbR$:
indeed, $\partial_\psi H=h\ge r$. The inverse of $H(r,\cdot)$ recovers the original data as
\[
 \psi=H^{-1}\bigl(r,H(r,\phi(r))+r^2u\bigr),\qquad
 \dot\psi=\frac{r^2\dot u}{h(r,\psi)}.
\]
\begin{prop}[Nonlinear transform]\label{prop:nonlinear-transform}
For every $\boldsymbol{\psi}\in\mathcal E_{0,n}$, we have
\begin{equation}\label{eq:data-energy-control}
 \norm{\boldsymbol{u}}_{\mathcal H}
       \lesssim 1+\sqrt{E(\boldsymbol{\psi})}.
\end{equation}
Moreover, the transformation \eqref{eq:transform} is a homeomorphism from
$\mathcal E_{0,n}$ onto $\mathcal H$.
\end{prop}

\begin{proof}
We first prove \eqref{eq:data-energy-control} and surjectivity,
then continuity in both directions. Injectivity follows from the
pointwise inverse above. For any
$\boldsymbol\psi,\boldsymbol\varphi\in\mathcal E_{0,n}$ with images
$\boldsymbol u,\boldsymbol v$, let
$\boldsymbol w=\boldsymbol u-\boldsymbol v$. By the definition \eqref{eq:transform}, we have
\begin{equation}\label{eq:transform-differences}
\begin{aligned}
 r^2\dot w&=h(r,\psi)\dot\psi-h(r,\varphi)\dot\varphi,\\
 \partial_r(r^2w)&=h(r,\psi)\psi_r-h(r,\varphi)\varphi_r
                  +H_r(r,\psi)-H_r(r,\varphi).
\end{aligned}
\end{equation}
Since $|H_{r\psi}|\le1$, Hardy's inequality implies
\[
 \int_0^\infty|H_r(r,\psi)-H_r(r,\varphi)|^2dr
 \le\int_0^\infty|\psi-\varphi|^2dr
 \le4\int_0^\infty r^2|\psi_r-\varphi_r|^2dr.
\]
Thus $\boldsymbol\psi,\boldsymbol\varphi\in\mathcal E_{0,n}$ implies $\partial_r(r^2w)\in L^2(0,\infty)$.
Also, $\psi,\varphi\to0$ as $r\downarrow0$, so $r^2w\to0$.
Integration by parts and Hardy's inequality therefore imply
\begin{equation}\label{eq:transform-energy-identity}
 \int_0^\infty \! r^4w_r^2dr 
 =\int_0^\infty \! |\partial_r(r^2w)|^2dr
       +2\int_0^\infty \! r^2w^2dr
 \le9\int_0^\infty|\partial_r(r^2w)|^2dr.
\end{equation}
Together with \eqref{eq:transform-differences}, this yields
\begin{align}
\|\boldsymbol w\|_{\mathcal H}
 \lesssim & 
 \left\|\frac{h(r,\psi)\psi_r-h(r,\varphi)\varphi_r}{r}\right\|_{L^2(\bbR^3)}
 +\left\|\frac{h(r,\psi)\dot\psi-h(r,\varphi)\dot\varphi}{r}\right\|_{L^2(\bbR^3)}\nonumber\\
 &+\|\psi_r-\varphi_r\|_{L^2(\bbR^3)}. \label{eq:forward-components}
\end{align}

Next, we claim that
\begin{equation}\label{eq:inverse-components}
\begin{aligned}
 &\left\|\frac{h(r,\psi)\psi_r-h(r,\varphi)\varphi_r}{r}\right\|_{L^2(\bbR^3)}
 +\left\|\frac{h(r,\psi)\dot\psi-h(r,\varphi)\dot\varphi}{r}\right\|_{L^2(\bbR^3)}\\
 &\qquad+\left\|\frac{\sin\psi-\sin\varphi}{r}\right\|_{L^2(\bbR^3)}
 +\left\|\frac{\sin^2\psi-\sin^2\varphi}{r^2}\right\|_{L^2(\bbR^3)}
 \lesssim\|\boldsymbol w\|_{\mathcal H}.
\end{aligned}
\end{equation}
By \eqref{eq:transform-energy-identity} and \eqref{eq:transform-differences}, \eqref{eq:inverse-components} it suffices to prove that 
\[ \left\|\frac{H_r(r,\psi)-H_r(r,\varphi)}{r}\right\|_{L^2(\bbR^3)}
 +\left\|\frac{\sin\psi-\sin\varphi}{r}\right\|_{L^2(\bbR^3)}
 +\left\|\frac{\sin^2\psi-\sin^2\varphi}{r^2}\right\|_{L^2(\bbR^3)} \]
is $O(\|\boldsymbol w\|_{\mathcal H})$. The bounds $|\cos\psi|\le1$, $|H_{r\psi}|\le1$, and $r\le h(r,\psi)$ imply
\[|\sin\psi-\sin\varphi| +|H_r(r,\psi)-H_r(r,\varphi)| \le 2|\psi-\varphi| \lesssim \left|\int_{\varphi}^{\psi} \frac{h(r,s)}{r}ds \right| = r |w|.\]
Moreover, $|\sin2\psi|\le\sqrt2 h(r,\psi)$ implies
\[\frac{|\sin^2\psi-\sin^2\varphi|}{r}\le \left| \int_{\varphi}^{\psi} \frac{|\sin 2s|}{r} ds \right| \lesssim \left|\int_{\varphi}^{\psi} \frac{h(r,s)}{r}ds \right|
 =r|w|.\]
The pointwise estimates above therefore
establish \eqref{eq:inverse-components}. We note that \eqref{eq:inverse-components} also holds for the pointwise
inverse images of arbitrary $\boldsymbol u,\boldsymbol v\in\mathcal H$. 

For $\boldsymbol\varphi=\boldsymbol\phi$, we have
$\boldsymbol v=0$. \eqref{eq:data-energy-control} follows
from \eqref{eq:forward-components}. Also,
\eqref{eq:inverse-components} shows that the inverse image of any
$\boldsymbol u\in\mathcal H$ has finite energy. For the endpoints, the radial Sobolev inequality implies
\[
\begin{aligned}
 \sqrt2\left|\int_0^\psi|\sin s|d s\right|
 &\le|H(r,\psi)|=r^2|u|
 \lesssim\|u\|_{\dot H^1(\bbR^5)}\sqrt r
 &&\text{for }0<r\le1,\\
 |\psi-n\pi|&\le r|u|
 \lesssim\|u\|_{\dot H^1(\bbR^5)}r^{-1/2}
 &&\text{for }r\ge2.
\end{aligned}
\]
Here we used $H(r,\psi)-H(r,n\pi)=H(r,\psi-n\pi)$ in the
second line. Thus $\psi(0)=0$ and $\psi(\infty)=n\pi$, which proves
surjectivity onto $\mathcal H$.

For forward continuity, we assume $d(\boldsymbol\psi_j,\boldsymbol\psi)\to0$ as $j\to \infty$, and let
$\boldsymbol u_j,\boldsymbol u$ be their transforms.
The radial Sobolev inequality in $\bbR^3$ implies
$\psi_j(r)\to\psi(r)$ for every $r>0$. We note that
\[
 \frac{h(r,\psi_j)\partial_r\psi_j}{r}
 =\frac r{h(r,\psi_j)}\partial_r\psi_j
  +\frac{2\sin\psi_j}{h(r,\psi_j)}
       \frac{\sin\psi_j\,\partial_r\psi_j}{r}.
\]
By the definition of the metric $d$, both $\partial_r\psi_j$ and
$\sin\psi_j\,\partial_r\psi_j/r$ converge in $L^2(\bbR^3)$
to the expressions with $\psi_j$ replaced by $\psi$.
The coefficients $r/h(r,\psi_j)$ and
$2\sin\psi_j/h(r,\psi_j)$ are uniformly bounded and converge pointwise.
Consequently, as $j\to \infty$,
\[
 \frac{h(r,\psi_j)\partial_r\psi_j}{r}
 \to\frac{h(r,\psi)\psi_r}{r}
 \quad\text{in }L^2(\bbR^3).
\]
The same computation with $\dot\psi_j$ in place of $\partial_r\psi_j$
proves convergence of $h(r,\psi_j)\dot\psi_j/r$. Hence, \eqref{eq:forward-components} with $\boldsymbol{\varphi}=\boldsymbol{\psi}_j$ implies $\boldsymbol u_j\to\boldsymbol u$ in $\mathcal H$.

Conversely, we assume $\boldsymbol u_j\to\boldsymbol u$ in $\mathcal H$,
and denote their inverse images by $\boldsymbol\psi_j,\boldsymbol\psi$. Again, the radial Sobolev
inequality with $|\psi_j-\psi|\le r|u_j-u|$ implies $\psi_j(r)\to\psi(r)$ for every $r>0$.
By \eqref{eq:inverse-components},
\[
 \frac{h(r,\psi_j)\partial_r\psi_j}{r}
 \to\frac{h(r,\psi)\psi_r}{r},\qquad
 \frac{h(r,\psi_j)\dot\psi_j}{r}
 \to\frac{h(r,\psi)\dot\psi}{r}
 \quad\text{in }L^2(\bbR^3).
\]
Since $r/h(r,\psi_j)$ and $\sin\psi_j/h(r,\psi_j)$
are uniformly bounded and converge pointwise, we obtain
convergence of $\partial_r\psi_j$, $\dot\psi_j$,
$\sin\psi_j\,\partial_r\psi_j/r$, and
$\sin\psi_j\,\dot\psi_j/r$ in $L^2(\bbR^3)$.
The remaining $\sin^2$ term in the metric $d$ is controlled by \eqref{eq:inverse-components}. Thus
$d(\boldsymbol\psi_j,\boldsymbol\psi)\to0$.
\end{proof}

We next compute the corresponding transformed equation.
The transformed nonlinearity has two useful growth bounds: the quadratic bound enters
the local theory, while the cubic bound supplies the exterior decay
needed later. The contribution from the fixed annulus is globally Lipschitz.
\begin{lem}[Nonlinear estimates]\label{lem:nonlinear-estimates}
For every smooth solution $\boldsymbol\psi=(\psi,\psi_t)$ of \eqref{eq:skyrme},
the transformed solution $\boldsymbol u=(u,u_t)$ defined by \eqref{eq:transform} satisfies
\begin{equation}\label{eq:relative-equation}
 u_{tt}-\Delta_5u=F(r,u),\quad \text{for }r>0,
\end{equation}
 where the nonlinearity $F$ is given by  
\begin{equation}\label{eq:full-source-formula}
 F(r,u)=\frac1{r^2}\biggl[2u
 +\frac{d^2}{ dr^2}H(r,\phi(r))-H_{rr}(r,\psi)
 -\frac{\sin2\psi}{h(r,\psi)}
               \left(1+\frac{\sin^2\psi}{r^2}\right)\biggr]
\end{equation}
with $\psi=H^{-1}(r,H(r,\phi(r))+r^2u)$. Moreover, $F$ satisfies:
\begin{itemize}
\item (Global derivative bound) Let $b\coloneqq \sin^2 \psi / r^2$. For $r>0$,
\begin{equation}\label{eq:vacuum-equation}
 F_u(r,u)=3b+\frac{b}{(1+2b)^2},\quad 0\le F_u(r,u)\le 4b \le \frac4{r^2}.
\end{equation}

\item (Bounds on $0<r\le1$ and $r\ge2$) We have $F(r,0)=0$, and
\begin{align}
 0\le F_u(r,u)&\le4\sqrt2|u|, &|F(r,u)|&\le2\sqrt2|u|^2,
 \label{eq:quadratic-source}\\
 0\le F_u(r,u)&\le4|u|^2, &|F(r,u)|&\le\tfrac43|u|^3.
 \label{eq:cubic-source}
\end{align}

\item (Bounds on $1<r<2$) We have
\[
 |F(r,0)|\lesssim 1,\qquad 0\le F_u(r,u)\le4,
\]
and consequently
\[
 |F(r,u)|\lesssim 1+|u|,\qquad
 |F(r,u)-F(r,v)|\le4|u-v|.
\]

\item (Global subcritical bounds) For every $r>0$, we have
\begin{align}
 |F(r,u)|&\lesssim|u|^2+\mathbf1_{1<r<2}(1+|u|),
                  \label{eq:global-source}\\
 |F(r,u)-F(r,v)|&\lesssim(|u|+|v|)|u-v|
                       +\mathbf1_{1<r<2}|u-v|.
                  \label{eq:source-lipschitz}
\end{align}
\end{itemize}
Near $r=0$, $F$ extends smoothly and evenly in $r$ for
bounded values of $u$.
\end{lem}
\begin{rem}\label{rem:WMYM interpolation}
Near $r=0$ or $r=\infty$, our nonlinear transformation \eqref{eq:transform} and the
transformed nonlinearity $F$ satisfy the following asymptotic equivalences
\[ u\sim
\begin{cases}
 \displaystyle\frac{\psi-\phi}{r},\\[6pt]
 \displaystyle\frac{\sqrt2\operatorname{sgn}(\psi)(1-\cos\psi)}{r^2},
\end{cases}
\;
F(r,u)\sim
\begin{cases}
 \displaystyle\frac{2ru-\sin(2ru)}{r^3},\\[6pt]
 \displaystyle\frac3{\sqrt2}u|u|-\frac{r^2}{2}u^3,
\end{cases}
\;
\begin{array}{@{}l@{}}
 \displaystyle\vphantom{\frac{2u}{r^2}-\frac{\sin(2ru)}{r^3}}
 \text{as }b\to0,\\[6pt]
 \displaystyle\vphantom{\frac3{\sqrt2}u|u|-\frac{r^2}{2}u^3}
 \text{as }b\to\infty.
\end{array}
\]
These asymptotics show that
the transformation \eqref{eq:transform} interpolates between supercritical wave map and subcritical Yang--Mills dynamics in \eqref{eq:WMYM 5D variables}.
\end{rem}
\begin{proof}
Here $\partial_r$ leaves the angle fixed, whereas $d/dr$ on
$H(r,\phi(r))$ denotes the total derivative. Substitution of
$H(r,\psi)=H(r,\phi(r))+r^2u$ into \eqref{eq:skyrme} cancels the derivative
nonlinearities, leaving \eqref{eq:relative-equation} with
the transformed nonlinearity \eqref{eq:full-source-formula}. At fixed $r$, differentiation in $u$, using $\psi_u=r^2/h(r,\psi)$,
yields \eqref{eq:vacuum-equation}; see also \cite[(5.10)]{Li21}.
The inequalities $0\le b/(1+2b)^2\le b$ and $b\le r^{-2}$ yield the remaining bounds in \eqref{eq:vacuum-equation}. 

For $0<r\le1$ or $r\ge2$, $\phi$ is respectively $0$ or
$n\pi$, and periodicity of $h$ implies
\[
 r^2u = H(r,\psi)-H(r,\phi)=H(r,\psi-\phi).
\]
At $u=0$, $H(r,\psi)=H(r,\phi)$, and thus $\psi=\phi$. Since $\phi',\phi'',\sin2\phi=0$, we have $F(r,0)=0$. Also,
\[
 |b_u|=\frac{|\sin 2\psi|}{r^2} \frac{r^2}{h(r,\psi)}=\frac{|\sin2\psi|}{h(r,\psi)}\le\sqrt2,
 \qquad b(r,0)=0.
\]
This proves $b\le\sqrt2|u|$. On the other hand, $\phi\equiv 0,n\pi$ and $h(r,\cdot)/r \ge 1$ imply
\[|\sin \psi|=|\sin (\psi-\phi)| \le |\psi-\phi| = \left| \int_{\phi}^{\psi} 1 ds \right| \le \left|\int_{\phi}^{\psi} \frac{h(r,s)}{r}ds \right| \le r |u|. \]
Hence, we have $b\le u^2$. The bounds $b\le\sqrt2|u|$ and $b\le u^2$, together with \eqref{eq:vacuum-equation} and integration
from $0$ to $u$, prove \eqref{eq:quadratic-source} and \eqref{eq:cubic-source}.

On $1<r<2$, boundedness of $F(r,0)$ and the estimate $F_u\le4$
from \eqref{eq:vacuum-equation} prove the stated annular bounds.
The annular bounds and \eqref{eq:quadratic-source} yield
\eqref{eq:global-source} and \eqref{eq:source-lipschitz}.

For smoothness at the origin, we use $\phi=0$ and set $g=\psi/r$. Then
\[
 u=\int_0^g\sqrt{1+2\bigl(\sin(rs)/r\bigr)^2}ds.
\]
The integral is smooth in $(r^2,g)$, with positive derivative in $g$. Its inverse $g=g(r,u)$ is therefore smooth in $(r^2,u)$.
Consequently, $b=(\sin(rg(r,u))/r)^2$ is smooth in $(r^2,u)$.
Since $F(r,0)=0$, integrating $F_u$ yields the extension.
\end{proof}

We note also that smooth compactly supported radial data $\boldsymbol u$
have smooth corotational inverse images, equal to $\boldsymbol\phi$
outside a compact set. By Proposition~\ref{prop:nonlinear-transform},
their inverse images form a dense subset of $\mathcal E_{0,n}$.

\subsection{Global well-posedness}
For $\boldsymbol u\in C(I;\mathcal H)$ with $u\in L^4_I L^4_x$
on a short interval $I$, the nonlinear estimates imply
$F(r,u)\in L^1_I L^2_x$, so we formulate the local problem as
\[
 \boldsymbol u(t)=S(t-t_0)\boldsymbol u(t_0)
 +\int_{t_0}^tS(t-s)(0,F(r,u(s)))d s
 \quad\text{for }t,t_0\in I.
\]

\begin{lem}[Local well-posedness]\label{lem:lwp}
Equation~\eqref{eq:relative-equation} is locally well-posed in
$\mathcal H$, with a lifespan depending only on the initial norm.
If $\|\boldsymbol{u}(t)\|_{\mathcal H}\le M$ on an interval $I$,
then for every $L<\infty$,
\begin{equation}\label{eq:strichartz}
 \sup_{[\tau,\tau+L]\subset I}\left(
 \norm{u}_{L^4_{[\tau,\tau+L]} L^4_x}
 +\norm{u}_{L^3_{[\tau,\tau+L]} L^{30/7}_x}
 +\norm{u}_{L^{7/3}_{[\tau,\tau+L]} L^{14/3}_x}\right)
 \lesssim_{M,L}1.
\end{equation}
\end{lem}

\begin{proof}
For $(q,p)=(4,4),(3,30/7),(7/3,14/3)$, the wave Strichartz inequality
\cite[Corollary~1.3]{KeelTao98} implies
\begin{equation}\label{eq:linear}
 \norm{\boldsymbol{u}}_{C(I;\mathcal H)}+\norm u_{L^q_I L^p_x}
 \lesssim\norm{\boldsymbol{u}(t_0)}_{\mathcal H}
               +\norm{F(r,u)}_{L^1_I L^2_x},\qquad\text{for }t_0\in I.
\end{equation}
Estimates \eqref{eq:global-source} and
\eqref{eq:source-lipschitz}, H\"older's inequality and the Sobolev
inequality on $1<r<2$ yield
\begin{align*}
 \norm{F(r,u)}_{L^1_I L^2_x}
 \lesssim & \, |I|^{1/2}\norm u_{L^4_I L^4_x}^2
   +|I|\bigl(1+\norm u_{L^\infty_I \dot H^1_x}\bigr),
 \\
 \norm{F(r,u)-F(r,v)}_{L^1_I L^2_x}
 \lesssim & \, |I|^{1/2}
   \bigl(\norm u_{L^4_I L^4_x}+\norm v_{L^4_I L^4_x}\bigr)
       \norm{u-v}_{L^4_I L^4_x} \\
       &+|I|\norm{u-v}_{L^\infty_I \dot H^1_x}.
\end{align*}
The standard contraction argument proves local well-posedness,
with a lifespan uniform on bounded subsets of $\mathcal H$.
Iteration of \eqref{eq:linear} proves \eqref{eq:strichartz}.
\end{proof}

The energy functional for \eqref{eq:relative-equation} is
\[
 \frac12\|\boldsymbol u\|_{\mathcal H}^2-\mathcal V(u),
 \qquad
 \mathcal V(u)=\int_0^\infty
       \left(\int_0^{u(r)}F(r,s)d s\right)r^4dr.
\]
Lemma~\ref{lem:nonlinear-estimates} and the
Sobolev inequality imply that $\mathcal V$ is well defined
and continuous on $\dot H^1_{\rm rad}(\bbR^5)$.

\begin{lem}[Energy identity]\label{lem:energy-identity}
For every $\boldsymbol\psi\in\mathcal E_{0,n}$ with transform
$\boldsymbol u$,
\begin{equation}\label{eq:transformed-hamiltonian}
 E(\boldsymbol\psi)
 =E(\boldsymbol\phi)+\frac12\|\boldsymbol u\|_{\mathcal H}^2
       -\mathcal V(u).
\end{equation}
\end{lem}

\begin{proof}
For smooth compactly supported radial data $\boldsymbol u$, let
\[
 \psi_s(r)
 =H^{-1}\bigl(r,H(r,\phi(r))+sr^2u(r)\bigr),
 \qquad\text{for }0\le s\le1.
\]
We set $\boldsymbol\psi_s=(\psi_s,0)$. Then
$\psi_0=\phi$, $\psi_1=\psi$, and $h(r,\psi_s)\partial_s\psi_s=r^2u$.
For a smooth compactly supported variation $\zeta$,
the first variation is
\[
 \left.\frac{d}{d\varepsilon}E(\psi+\varepsilon\zeta,0)
       \right|_{\varepsilon=0}
 =-\int_0^\infty h(r,\psi)
       [\Delta_5u+F(r,u)]\zeta\,r^2dr.
\]
Hence, we obtain
\begin{align*}
      \frac{d}{ds}E(\boldsymbol\psi_s) & = -\int_0^{\infty}  h(r,\psi_s) \bigl[\Delta_5 (su) +F(r,su)\bigr] \partial_s\psi_s r^2 dr \\
      &=\!\int_0^\infty
       \!\bigl[-\Delta_5(su)-F(r,su)\bigr]u\,r^4dr=\frac{d}{ds}\!\left[
       \frac{s^2}{2}\int_0^\infty u_r^2r^4dr-\mathcal V(su)
       \right].
\end{align*}
Integration from $s=0$ to $s=1$ with $\mathcal V(0)=0$ yields
the identity \eqref{eq:transformed-hamiltonian} for $(\psi,0)$ and $(u,0)$.
The kinetic energies agree by \eqref{eq:transform}, so the identity holds for
smooth compactly supported $\boldsymbol u$.
The identity extends to all $\boldsymbol u\in\mathcal H$
by density, Proposition~\ref{prop:nonlinear-transform},
and continuity of $E$ and $\mathcal V$.
\end{proof}

\begin{prop}[Global well-posedness]\label{prop:transformed-gwp}
For every initial datum $\boldsymbol u_0\in\mathcal H$, the maximal
solution $\boldsymbol u=(u,u_t)$ given by Lemma~\ref{lem:lwp}
conserves energy:
\[
 \frac12\|\boldsymbol u(t)\|_{\mathcal H}^2-\mathcal V(u(t))
 =\frac12\|\boldsymbol u_0\|_{\mathcal H}^2-\mathcal V(u_0)
 \qquad\text{for }t\in I.
\]
Consequently, $I=\bbR$ and
\begin{equation}\label{eq:uniform-energy}
 \sup_{t\in\bbR}\norm{\boldsymbol u(t)}_{\mathcal H}
 \lesssim_E1,\qquad E=E(\boldsymbol\psi_0).
\end{equation}
Equation~\eqref{eq:relative-equation} is therefore globally
well-posed in $\mathcal H$.
Moreover, uniqueness holds among all distributional solutions
in $C(I;\mathcal H)$ on any time interval $I$ containing zero.
\end{prop}
\begin{proof}
For smooth solutions, \eqref{eq:relative-equation} and
\eqref{eq:transformed-hamiltonian} imply energy conservation.
The general case follows by smooth approximation, using
Proposition~\ref{prop:persistence}, continuous dependence,
and continuity of $\mathcal V$.
The bound \eqref{eq:data-energy-control} and the uniform local
lifespan then imply global existence.

For uniqueness, let $\boldsymbol u=(u,u_t),\boldsymbol v=(v,v_t)\in C(I;\mathcal H)$ be distributional solutions with $\boldsymbol u(t_0)=\boldsymbol v(t_0)$. The first component $w$ of $\boldsymbol w=\boldsymbol u-\boldsymbol v$ lies in $C(I;H^1)$, since
$w(t)=\int_{t_0}^t w_t(s)d s$. For a bounded interval $J\subset I$ containing $t_0$,
the Duhamel formula, the wave Strichartz inequality, and \eqref{eq:source-lipschitz} yield
\begin{align*}
 \norm w_{L^4_J L^{20/7}_x}
 & \lesssim\norm{F(r,u)-F(r,v)}_{L^{4/3}_J L^{20/13}_x} \\
 & \lesssim |J|^{1/2}
 \left(1+\norm u_{L^\infty_J\dot H^1_x}
          +\norm v_{L^\infty_J\dot H^1_x}\right)
 \norm w_{L^4_J L^{20/7}_x}.
\end{align*}
By $H^1\subset L^{20/7}$, $\norm w_{L^4_JL^{20/7}_x}$ is finite. We absorb the right-hand side
on a short interval and repeat the argument to prove uniqueness
on $I$.
\end{proof}

We now prove Theorem~\ref{thm:gwp} in the original variables.
\begin{proof}[Proof of Theorem~\ref{thm:gwp}]
Let $\boldsymbol\psi_0=(\psi_0,\dot\psi_0)\in\mathcal E_{0,n}$.
By Proposition~\ref{prop:nonlinear-transform}, the transformed datum belongs
to $\mathcal H$ and hence determines a unique global solution $\boldsymbol u$
by Proposition~\ref{prop:transformed-gwp}.
The inverse transformation defines
$\boldsymbol\psi\in C(\bbR;\mathcal E_{0,n})$
with initial data $\boldsymbol\psi_0$.
We first prove that $\boldsymbol\psi$ satisfies \eqref{eq:skyrme} in the sense of distributions
\begin{equation}\label{eq:energy-weak-form}
 \begin{aligned}
 \iint \! \!\big[
 &h^2\psi_t\xi_t-h^2\psi_r\xi_r
 +\sin(2\psi)(\psi_t^2-\psi_r^2)\xi-\sin(2\psi)(1+\tfrac{\sin^2\psi}{r^2})\xi
 \big]drdt=0,
 \end{aligned}
\end{equation}
where $\xi\in C_c^\infty(\bbR\times[0,\infty))$ and $\xi/r$
extends smoothly and evenly across $r=0$. The chain rule for
$H(r,\psi)=H(r,\phi)+r^2u$ yields
\[
 \psi_t=\frac{r^2u_t}{h(r,\psi)}=\dot\psi
 \qquad\text{for }r>0.
\]
For $\xi$ supported away from the origin, we take
$\varphi=h(r,\psi)\xi/r^2$ as a test function in the weak form
of \eqref{eq:relative-equation}. We approximate $\varphi$ by smooth test functions,
both uniformly and in $H^1$ with respect to $(t,r)$. This justifies the choice of $\varphi$ as a test function.
The chain rule and the formula for $F$ then yield
\eqref{eq:energy-weak-form}. For a general test function $\xi$, we fix a compact time interval $I$ containing its time support. We then apply \eqref{eq:energy-weak-form}
to $(1-\chi_\varepsilon)\xi$. Since $\xi=O(r)$ near the origin, 
the only term containing $\chi_\varepsilon'$ satisfies
\[
 \left|\iint h^2\psi_r\xi(\chi_\varepsilon)' drdt\right|
 \lesssim_{\xi,I}\varepsilon^{1/2}
 \left\|\frac{h\psi_r}{r}\right\|_{L^\infty_I L^2(\bbR^3)} \lesssim \varepsilon^{1/2} \sup_{t\in I}E(\boldsymbol\psi(t))^{1/2}.
\]
All other terms are integrable since $\xi=O(r)$ and
$\sup_{t\in I}E(\boldsymbol\psi(t))<\infty$.
By the dominated convergence theorem, we may pass to the limit
$\varepsilon\to0$ and obtain \eqref{eq:energy-weak-form}.

Now, we prove the uniqueness part. Let $I$ be a compact time interval and
$\boldsymbol\psi=(\psi,\psi_t)\in C(I;\mathcal E_{0,n})$
satisfy \eqref{eq:energy-weak-form}.
Proposition~\ref{prop:nonlinear-transform} and the chain rule
imply $\boldsymbol u=(u,u_t)\in C(I;\mathcal H)$.
The same argument, now with
$\xi=r^2\varphi/h(r,\psi)$, proves
\eqref{eq:relative-equation} for $r>0$. To include $r=0$, we test with
$(1-\chi_{\varepsilon})\varphi$.
The only term containing $\chi_\varepsilon'$ satisfies
\[
 \left|\int_I\int_0^\infty
       \varphi u_r\chi_\varepsilon' r^4drdt\right|
 \lesssim_{\varphi,I}\varepsilon^{3/2}
 \|u_r\|_{L^\infty_I L^2(\bbR^5)} \lesssim \varepsilon^{3/2} \sup_{t\in I}\norm{\boldsymbol u(t)}_{\mathcal H}.
\]
The other terms converge by the dominated convergence theorem; the quadratic bound on $F$ implies
$F(r,u)\in L^1_{\mathrm{loc}}(I\times\bbR^5)$.
Thus \eqref{eq:relative-equation} holds on $\bbR^5$,
and Proposition~\ref{prop:transformed-gwp} proves uniqueness.
Energy conservation and continuous dependence follow from
Propositions~\ref{prop:nonlinear-transform}
and~\ref{prop:transformed-gwp}, together with
\eqref{eq:transformed-hamiltonian}.
\end{proof}

 Let $\boldsymbol Q=(Q_{0,n},0)$ be the stationary solution in $\mathcal E_{0,n}$, with first component $Q=Q_{0,n}$. We set
\begin{equation}\label{eq:transformed-stationary-profile}
 \boldsymbol{u}_Q=(u_Q,0),\qquad u_Q(r)=\frac{H(r,Q(r))-H(r,\phi(r))}{r^2}.
\end{equation}
By Proposition~\ref{prop:nonlinear-transform} and
Theorem~\ref{thm:stationary}, $\boldsymbol u_Q$ is the unique
stationary solution of \eqref{eq:relative-equation} in $\mathcal H$. Now, we state soliton resolution for $\boldsymbol u$ in $\mathcal H$.
\begin{prop}[Soliton resolution]\label{prop:transformed-resolution}
For every global solution $\boldsymbol u$ in Proposition~\ref{prop:transformed-gwp},
there is a radial free wave $\boldsymbol u_L=(u_L,\partial_tu_L)\in C(\bbR;\mathcal H)$ on $\bbR^{1+5}$ such that
\begin{equation}\label{eq:transformed-resolution}
 \boldsymbol{u}(t)=\boldsymbol{u}_Q+\boldsymbol{u}_L(t)+o_{\mathcal H}(1)
 \qquad\text{as }t\to+\infty.
\end{equation}
\end{prop}

We prove Proposition~\ref{prop:transformed-resolution}, and then prove
Theorem~\ref{thm:main} in Section~\ref{sec:transformed-completion}.
\section{Compactness and radiation}\label{sec:radiation}

We prove in this section that, after subtraction of a free wave, the solution has a precompact forward trajectory in $\mathcal H$.

\begin{prop}[Compactness]\label{prop:compactness}
Let $\boldsymbol u$ be a solution in Proposition~\ref{prop:transformed-gwp}.
There is a radial finite-energy free wave $\boldsymbol u_L$ such that
\[
 \{\boldsymbol u(t)-\boldsymbol u_L(t):t\ge0\}
 \quad\text{is precompact in }\mathcal H.
\]
\end{prop}

We first establish local precompactness in $\mathcal H_B$, then
construct a free wave and obtain exterior tightness of the remainder.
We then prove Proposition~\ref{prop:compactness} at the end
of this section. Throughout this section, $\boldsymbol u$ is fixed, and constants
may depend on the conserved energy $E$.

\subsection{Local compactness}
For $R>2$, the radial Sobolev inequality and the cubic bound
\eqref{eq:cubic-source} imply
\begin{equation}\label{eq:exterior-source}
 \norm{\mathbf1_{r>R}F(r,u(t))}_{L^2(\bbR^5)}
 \lesssim_E R^{-2}.
\end{equation}
For the local norm in \eqref{eq:local-norm}, Hardy's inequality
implies $\|\boldsymbol f\|_{\mathcal H_B}
\lesssim\|\boldsymbol f\|_{\mathcal H}$ uniformly for $B>0$.

\begin{prop}[Local compactness]\label{prop:local-compactness}
For every $B>0$, $\{\boldsymbol{u}(t):t\geq0\}$ is precompact in $\mathcal H_B$.
\end{prop}

\begin{proof}
By continuity of $t\mapsto\boldsymbol u(t)$
in $\mathcal H$, the orbit on each bounded time interval is compact.
Let $B>0$. For $L>B+2$ and $t>L$, the Duhamel formula yields 
\[
 \boldsymbol u(t)=S(t)\boldsymbol u(0)
  +\!\int_0^{t-L}\!S(t-s)(0,F(r,u(s)))d s
  +\!\int_{t-L}^t \!S(t-s)(0,F(r,u(s)))d s.
\]
By the strong Huygens principle, the linear propagator term tends to zero in $\mathcal H_B$ as $t\to \infty$. 

For the second term, we set $\sigma=t-s\ge L$.
By the strong Huygens principle, only the values of
$F(r,u(t-\sigma))$ on $\sigma-B<r<\sigma+B$
contribute to $S(\sigma)(0,F(r,u(t-\sigma)))$ in $r<B$. Thus, \eqref{eq:exterior-source} yields
\begin{align*}
\left\|\int_0^{t-L}S(t-s)(0,F(r,u(s)))d s\right\|_{\mathcal H_B}
 &\lesssim\int_L^t\|\mathbf1_{r>\sigma-B}F(r,u(t-\sigma))\|_{L^2(\bbR^5)} d\sigma
 \\
 &\lesssim_{E} \frac{1}{L-B}.
\end{align*}
For sufficiently large $L$ and then sufficiently large $t$,
the first two terms are arbitrarily small in $\mathcal H_B$.
For each fixed $L>B+2$, let $\tau=t-L$ and $\sigma=s-\tau$. We claim precompactness of the forward trajectory of the last integral, 
\begin{equation*}
\begin{aligned}
 \left\{\int_0^L S(L-\sigma)(0,F(r,u(\tau+\sigma)))d\sigma:
                         \tau>0\right\}.
\end{aligned}
\end{equation*}
The Duhamel integral defines a bounded linear map from
$L^1_{[0,L]}L^2_x$ to $\mathcal H$. Precompactness of the integral
family therefore follows from precompactness of
$\{F(r,u(\tau+\cdot,\cdot)):\tau\ge0\}$ in $L^1_{[0,L]}L^2_x$.
We divide the spatial domain into three regions. Let $0<a<1$ and $R>2$.

(i) $r<a$: The quadratic bound on $F$, H\"older's inequality, and Lemma~\ref{lem:lwp} imply
\begin{align*}
\|\mathbf1_{r<a}F(r,u(\tau+\cdot))\|_{L^1_{[0,L]}L^2_x}
 &\lesssim\int_0^L\|\mathbf1_{r<a}\|_{L_x^{30}}
                  \|u(\tau+s)^2\|_{L_x^{15/7}}d s\\
 &\lesssim a^{1/6}L^{1/3}
       \|u\|_{L^3_{[\tau,\tau+L]} L^{30/7}_x}^2
 \lesssim_{E,L}a^{1/6}.
\end{align*}

(ii) $a\le r \le R$: By \eqref{eq:uniform-energy}, the family
$\{u(\tau+\cdot,\cdot):\tau\ge0\}$ is bounded in
\[
 L^\infty([0,L];H^1[a,R])
 \cap W^{1,\infty}([0,L];L^2[a,R]).
\]
Sobolev and Gagliardo--Nirenberg inequalities imply
equicontinuity in $(s,r)$. This family is therefore precompact
in $C([0,L]\times[a,R])$ by the Arzel\`a--Ascoli theorem.
By continuity of $F$, the set
$\{\mathbf1_{a\le r \le R}F(r,u(\tau+\cdot,r)):\tau\ge0\}$
is precompact in $L^1_{[0,L]}L^2_x$.

(iii) $r>R$: \eqref{eq:exterior-source} directly yields
\[
 \|\mathbf1_{r>R}F(r,u(\tau+\cdot))\|_{L^1_{[0,L]} L^2_x}
 \lesssim_{E}LR^{-2}.
\]

(i)--(iii) prove precompactness in
$L^1_{[0,L]}L^2_x$, which completes the proof of Proposition~\ref{prop:local-compactness}.
\end{proof}

\subsection{Exterior radiation}

We next construct a free wave $\boldsymbol u_L$ and control the exterior
energy of $\boldsymbol u-\boldsymbol u_L$ at large times.
We recall the linear radiation representation and channel inequality.

\begin{lem}[Free radial waves, \cite{DKM19Radiation,KLS14}]\label{lem:channel}
Every radial finite-energy free wave $\boldsymbol W$ has a unique radiation profile
$G_W\in L^2(\bbR)$ such that
\[
 \left\|W_t-\frac{G_W(r-t)}{r^2}\right\|_{L^2(\bbR^5)}
 +\left\|W_r+\frac{G_W(r-t)}{r^2}\right\|_{L^2(\bbR^5)}
 \to0\qquad\text{as }t\to+\infty.
\]
The map $\boldsymbol W(0)\mapsto G_W$ is onto, and
\[
 \norm{\boldsymbol W(0)}_{\mathcal H}^2=2\norm{G_W}_{L^2(\bbR)}^2.
\]
For $R>0$, let $\pi_R$ denote the orthogonal projection in
$\mathcal H(r>R)$ onto
\[
 P(R)=\operatorname{span}\{(r^{-3},0),(0,r^{-3})\},
\]
and let $\pi_R^\perp=I-\pi_R$.
Then
\begin{equation}\label{eq:free-channel}
 \max_{\sigma\in\{-,+\}}\lim_{t\to\sigma\infty}
   \norm{\boldsymbol{W}(t)}_{\mc H(r>R+|t|)}^2
       \gtrsim \norm{\pi_R^\perp\boldsymbol{W}(0)}_{\mc H(r>R)}^2.
\end{equation}
\end{lem}

The radiation representation, surjectivity, and norm identity follow by
restriction to radial waves from \cite[Theorem~2.1]{DKM19Radiation},
whose proof is in Appendices~A and~B of that paper. The channel inequality \eqref{eq:free-channel} follows from
\cite[Proposition~4.1]{KLS14}.

\medskip

We record the projection formulas needed for the exterior channel
argument, as in \cite[(5.12), (5.13)]{LawrieRodriguez19}.
For $\boldsymbol f=(f,\dot f)\in\mathcal H(r>R)$, the projection coefficients onto $P(R)$ are obtained from
\begin{equation}\label{eq:projection-coefficients}
 v_0(r)=r^3f(r),\qquad v_1(r)=r\int_r^\infty s\dot f(s)d s.
\end{equation}
The projection and the corresponding norm decomposition are
\begin{equation}\label{eq:projection-split}
 \pi_R\boldsymbol f=\left( \frac{v_0(R)}{r^{3}},\frac{v_1(R)}{r^{3}}\right),\quad
 \|\boldsymbol f\|_{\mathcal H(r>R)}^2
 =\|\pi_R^\perp\boldsymbol f\|_{\mathcal H(r>R)}^2+\frac{3v_0(R)^2}{R^3}+\frac{v_1(R)^2}{R}.
\end{equation}
Integration by parts yields
\[
 \|\pi_R^\perp\boldsymbol f\|_{\mathcal H(r>R)}^2=\int_R^\infty
       \left(\frac{|\partial_rv_0|^2}{r^2}+|\partial_rv_1|^2\right)d r.
\]
By the Cauchy--Schwarz inequality on $[R,r]$, we obtain
\begin{equation}\label{eq:dyadic-increments}
 \frac{|v_0(r)-v_0(R)|}{R^{3/2}}+\frac{|v_1(r)-v_1(R)|}{R^{1/2}}
 \lesssim\|\pi_R^\perp\boldsymbol f\|_{\mathcal H(r>R)},\quad\text{for }R\le r\le2R.
\end{equation}

\begin{prop}[Exterior radiation]\label{prop:radiation}
There is a unique radial finite-energy free wave $\boldsymbol u_L$ such that
\begin{equation}\label{eq:tail-tightness}
 \lim_{R\to\infty}\limsup_{t\to\infty}
 \norm{\boldsymbol u(t)-\boldsymbol u_L(t)}_{\mc H(r>R)}=0.
\end{equation}
\end{prop}

\begin{proof}
\underline{Step 1}: Radiation.
We construct a free wave $\boldsymbol u_L$ such that,
for every fixed $A\in\bbR$,
\begin{equation}\label{eq:cone-radiation}
 \norm{\boldsymbol{u}(t)-\boldsymbol{u}_L(t)}_{\mc H(r>t+A)}\to0\qquad\text{as }t\to+\infty.
\end{equation}
We follow \cite[Proposition~3.4, Steps~2 and~3]{DuyckaertsYang21}.
Let $A\in\bbR$ and let $T\ge0$ satisfy $T+A>2$.
For $t\ge T$, let $\boldsymbol u^A=(u^A,\partial_tu^A)$ solve
\[
 u^A_{tt}-\Delta_5u^A=\mathbf1_{r>t+A}F(r,u(t)),\qquad
 \boldsymbol u^A(T)=\boldsymbol u(T).
\]
Finite speed of propagation implies
$\boldsymbol u^A(t,r)=\boldsymbol u(t,r)$ for $t\ge T$ and $r>t+A$.
By \eqref{eq:exterior-source},
\[
 \int_T^\infty\|\mathbf1_{r>t+A}F(r,u(t))\|_2\,dt
 \lesssim_E\int_T^\infty(t+A)^{-2}\,dt
 =(T+A)^{-1}.
\]
Thus $\boldsymbol u^A$ scatters in $\mathcal H$ to a radial free
wave $\boldsymbol u_L^A$, with
$\|\boldsymbol u_L^A(0)\|_{\mathcal H}\lesssim_E1$
uniformly in $A$.

For $A'>A$, both free waves approximate $\boldsymbol u(t)$ on
$r>t+A'$. Lemma~\ref{lem:channel} therefore implies that their
radiation profiles satisfy
\[
 G^A=G^{A'}\quad\text{on }(A',\infty),\qquad
 \|G^A\|_{L^2(\bbR)}\lesssim_E1.
\]
The compatible profiles define $G\in L^2(\bbR)$ satisfying
$G|_{(A,\infty)}=G^A|_{(A,\infty)}$ for every $A\in\bbR$.
The free wave $\boldsymbol u_L$ with profile $G$ then satisfies
\eqref{eq:cone-radiation} by Lemma~\ref{lem:channel}.

\medskip
\underline{Step 2}: Tightness.
For $\boldsymbol w=\boldsymbol u-\boldsymbol u_L$, it suffices to show
\begin{equation}\label{eq:quantitative-tightness}
 \limsup_{T\to\infty}\norm{\boldsymbol w(T)}_{\mathcal H(r>R)}^2
 \lesssim_E R^{-1},\qquad\text{for }R\ge4.
\end{equation}
For $T\ge0$ and $R>2$, let $\boldsymbol W(s)=S(s)\boldsymbol w(T)$.
Since $w_{tt}-\Delta_5w=F(r,u)$, the Duhamel formula and finite
speed of propagation, together with \eqref{eq:linear} and
\eqref{eq:exterior-source}, imply
\begin{align}
 \norm{\boldsymbol W(s)-\boldsymbol w(T+s)}_{\mathcal H(r>R+|s|)}
 &\lesssim \left|\int_{0}^s
 \|\mathbf1_{r>R+|\tau|}F(r,u(T+\tau))\|_2 d\tau \right| \lesssim_E \frac{1}{R},\label{eq:exterior comparison}
\end{align}
for $-T\le s<\infty$. \eqref{eq:cone-radiation} with $t=T+s$ and $A=R-T$ implies
\[
 \lim_{s\to+\infty}\|\boldsymbol W(s)\|_{\mathcal H(r>R+s)}
 \lesssim_E R^{-1}.
\]
For the case $s\to -\infty$, monotonicity of the exterior
free-wave energy and \eqref{eq:exterior comparison} at $s=-T$ imply
\[
 \lim_{s\to-\infty}\|\boldsymbol W(s)\|_{\mathcal H(r>R+|s|)}
 \le\|\boldsymbol W(-T)\|_{\mathcal H(r>R+T)}
 \lesssim_E\|\boldsymbol w(0)\|_{\mathcal H(r>R+T)}+R^{-1}.
\]
Lemma~\ref{lem:channel} therefore implies
\begin{equation}\label{eq:remainder-defect}
 \norm{\pi_R^\perp\boldsymbol w(T)}_{\mathcal H(r>R)}
 \lesssim_E R^{-1}+\norm{\boldsymbol w(0)}_{\mathcal H(r>R+T)}.
\end{equation}

Let $v_0(T,r),v_1(T,r)$ be the coefficients
\eqref{eq:projection-coefficients} of $\boldsymbol w(T)$. For fixed $R$, we have $\|\boldsymbol w(0)\|_{\mathcal H(r>R+T)}=o_{T\to \infty}(1)$. \eqref{eq:remainder-defect} and
dyadic summation of \eqref{eq:dyadic-increments} from $R_0=4$,
with $\|\boldsymbol w(T)\|_{\mathcal H}\lesssim_E1$, yield
\[
 \limsup_{T\to\infty}|v_0(T,R)|\lesssim_E1+\sqrt R,\quad
 \limsup_{T\to\infty}|v_1(T,R)|\lesssim_E1,\quad \text{for }R\ge4.
\] The decomposition
\eqref{eq:projection-split} and \eqref{eq:remainder-defect} now yield \eqref{eq:quantitative-tightness}:
\[
 \limsup_{T\to\infty}\norm{\boldsymbol w(T)}_{\mathcal H(r>R)}^2
 \lesssim_E \frac{1}{R^{2}}+\frac{1+R}{R^3}+\frac{1}{R}
 \lesssim_E \frac{1}{R},\quad\text{for }R\ge4.
\]
Uniqueness follows from free wave energy conservation,
the strong Huygens principle, and \eqref{eq:tail-tightness}.
\end{proof}
\begin{proof}[Proof of Proposition~\ref{prop:compactness}]
Let $\boldsymbol u_L$ be the free wave from Proposition~\ref{prop:radiation}.
Proposition~\ref{prop:local-compactness} and the strong Huygens principle imply that the forward trajectory of $\boldsymbol w=\boldsymbol u-\boldsymbol u_L$
is precompact in every $\mathcal H_B$.
For every sequence $t_j\to\infty$, the sequence $\boldsymbol w(t_j)$
has a subsequence converging on each fixed ball. \eqref{eq:tail-tightness} makes this
subsequence Cauchy in $\mathcal H$. On finite
time intervals, compactness follows from continuity of
$t\mapsto\boldsymbol w(t)$ in $\mathcal H$, and the proof is complete.
\end{proof}

\section{Rigidity}\label{sec:rigidity}

In this section, we classify solutions with precompact trajectories. We then deduce Proposition~\ref{prop:transformed-resolution}
and Theorem~\ref{thm:main}.

\begin{prop}\label{prop:rigidity}
Let $\boldsymbol u_*=(u_*,\partial_tu_*)\in C(\bbR;\mathcal H)$ be a solution
to \eqref{eq:relative-equation}. If
\[
 \{\boldsymbol u_*(t):t\in\bbR\} \quad \text{is precompact in }\mathcal{H},
\]
then $\boldsymbol u_*(t)=\boldsymbol u_Q$ for every $t\in\bbR$, where
$\boldsymbol u_Q$ is defined in \eqref{eq:transformed-stationary-profile}.
\end{prop}

\begin{proof}[Sketch of proof]
We follow the rigidity argument of
\cite[Sections~5.3--5.5]{LawrieRodriguez19}. Precompactness directly implies
\begin{equation}\label{eq:hypotheses}
 \sup_{t\in\bbR}\norm{\boldsymbol u_*(t)}_{\mathcal H}\lesssim 1,
 \qquad
 \lim_{R\to\infty}\sup_{t\in\bbR}
       \norm{\boldsymbol u_*(t)}_{\mathcal H(r>R)}=0.
\end{equation}

\medskip\noindent
\underline{Step 1}: Channel estimate.
We prove, uniformly in $t$ for sufficiently large $R$,
\begin{equation}\label{eq:channel-estimate}
 \|\pi_R^\perp\boldsymbol u_*(t)\|_{\mathcal H(r>R)}
 \lesssim R^{-1}\|\boldsymbol u_*(t)\|_{\mathcal H(r>R)}^3.
\end{equation}
For $t\in\bbR$, let $\boldsymbol W(s)=S(s)\boldsymbol u_*(t)$. By \eqref{eq:free-channel} and \eqref{eq:hypotheses}, we obtain
\begin{align}
      \|\pi_R^\perp\boldsymbol u_*(t)\|_{\mathcal H(r>R)} &\lesssim \max_{\sigma\in\{-,+\}}\lim_{s\to\sigma\infty}
   \norm{S(s)\boldsymbol u_*(t)}_{\mc H(r>R+|s|)} \nonumber \\
   &\le  \sup_s\|\boldsymbol u_*(t+s)-S(s)\boldsymbol u_*(t)\|_{\mathcal H(r>R+|s|)}.\label{eq:orthogonal difference}\end{align}
The Duhamel formula with the cubic bound \eqref{eq:cubic-source} and the radial Sobolev inequality implies 
\begin{align*}
\eqref{eq:orthogonal difference} &\lesssim \max_{\sigma\in\{-,+\}}\int_{0}^{\infty} \|\mathbf1_{r>R+|\tau|}F(r,u_*(t+\sigma\tau))\|_2 d\tau \\
& \lesssim R^{-1}{\sup_s\|\boldsymbol u_*(t+s)\|_{\mathcal H(r>R+|s|)}^3 }.
\end{align*}
By the triangle inequality and the above estimate, we have
\[
 \sup_s\|\boldsymbol u_*(t+s)\|_{\mathcal H(r>R+|s|)}
 \lesssim \|\boldsymbol u_*(t)\|_{\mathcal H(r>R)}
 +R^{-1}\sup_s\|\boldsymbol u_*(t+s)\|_{\mathcal H(r>R+|s|)}^3.
\]
Hence, \eqref{eq:hypotheses} implies \eqref{eq:channel-estimate} for sufficiently large $R$.

\medskip\noindent
\underline{Step 2}: Projection coefficients. We claim that the coefficients \eqref{eq:projection-coefficients} of $\boldsymbol u_*(t)$ satisfy
\begin{equation}\label{eq:tail-limits}
 v_0(t,r)=\ell_0(t)+O(r^{-1}),\qquad v_1(t,r)=O(r^{-2}), \quad \text{as }r\to \infty,
\end{equation}
uniformly in time, with $\ell_0$ bounded.
By \eqref{eq:channel-estimate} and the uniform energy bound \eqref{eq:uniform-energy},
$\|\pi_R^\perp\boldsymbol u_*(t)\|_{\mathcal H(r>R)}=O(R^{-1})$.
Dyadic summation of \eqref{eq:dyadic-increments} from a fixed radius
implies $v_0(t,r)=O(r^{1/2})$ and $v_1(t,r)=O(1)$.
The decomposition \eqref{eq:projection-split} and \eqref{eq:channel-estimate} refine estimates of $\boldsymbol{u}_*$ as follows:
\[
 \|\boldsymbol u_*(t)\|_{\mathcal H(r>R)}=O(R^{-1/2}),\qquad
 \|\pi_R^\perp\boldsymbol u_*(t)\|_{\mathcal H(r>R)}=O(R^{-5/2}).
\]
Again, dyadic summation proves \eqref{eq:tail-limits} for $v_0$
and $v_1(t,r)=\ell_1(t)+O(r^{-2})$, uniformly in time, with bounded limits.
In particular, $u_*(t,r)=O(r^{-3})$, and so $F(r,u_*)=O(r^{-9})$. Together with the definition of $v_1$ in \eqref{eq:projection-coefficients} and averaging on $[R,2R]$, as in the proof of \cite[Claim 5.13]{LawrieRodriguez19}, we obtain
\[
 |\ell_1(t_2)-\ell_1(t_1)| \lesssim R^{-2} (1+|t_2-t_1|).
\]
As $R\to \infty$, we get $\ell_1(t)\equiv \ell_1$. Again, by the definition of $v_1$ and $u_*(t,r)=O(r^{-3})$, we have
\[
 |\ell_1|
 \le \frac RT\int_R^\infty |
 u_*(T,r)-u_*(0,r)|rdr+O(R^{-2}) \lesssim T^{-1}+R^{-2}.
\]
Taking $T = R$ and $R\to \infty$, we prove \eqref{eq:tail-limits}.

\medskip\noindent
\underline{Step 3}: Compact support.
For $\ell=\ell_0(0)$, we construct a stationary profile $\boldsymbol{u}_{\ell}=(u_\ell,0)$ with
$r^3u_\ell(r)\to\ell$ and prove that $\boldsymbol u_*(0)= \boldsymbol{u}_{\ell}$
on $r>R$ for large $R$. By a standard contraction argument with \eqref{eq:cubic-source}, for sufficiently large $R_1>2$, one can construct $u_\ell$ satisfying
\begin{equation}\label{eq:stationary-tail}
 u_\ell(r)=\frac{\ell}{r^3}
 -\int_r^\infty\frac1{s^4}
       \int_s^\infty \rho^4F(\rho,u_\ell(\rho))\,d\rho ds,
 \qquad r\ge R_1.
\end{equation}
Since $|F_u(r,\cdot )|\le4r^{-2}$, $u_\ell$ is defined for
all $r>0$. By \eqref{eq:vacuum-equation}, \eqref{eq:cubic-source}, $\boldsymbol{\widetilde u}_* \coloneqq \boldsymbol u_*-\boldsymbol{u}_{\ell}$ satisfies
\begin{align}
 \partial_t^2\widetilde u_*-\Delta_5\widetilde u_*=V(t,r)\widetilde u_*,\qquad
 &V(t,r)=\int_0^1F_u(r,u_\ell(r)+s\widetilde u_*(t,r))d s, \nonumber \\
|V(t,r)|\le4r^{-2},\qquad &|V(t,r)|\lesssim r^{-6}\quad\text{for }r\ge R_1.\label{eq:exterior-potential}
\end{align}
Precompactness of $\boldsymbol{u}_*$ and the tail decay of $u_{\ell}$ imply, for every $R>0$,
\begin{equation}\label{eq:difference-nonradiative}
 \lim_{t\to\pm\infty}\|\boldsymbol{\widetilde u}_*(t)\|_{\mathcal H(r>R+|t|)}=0.
\end{equation}
Repeating \underline{Step~1} with initial data
$\boldsymbol{\widetilde u}_*(0)$ and forcing $V\widetilde u_*$
yields
\begin{equation}\label{eq:second-channel}
 \|\pi_R^\perp\boldsymbol{\widetilde u}_*(0)\|_{\mathcal H(r>R)}
 \lesssim R^{-4}\|\boldsymbol{\widetilde u}_*(0)\|_{\mathcal H(r>R)}.
\end{equation}
Let $\widetilde v_0,\widetilde v_1$ be the projection coefficients of
$\boldsymbol{\widetilde u}_*(0)$. By \eqref{eq:tail-limits} and
\eqref{eq:stationary-tail}, $\widetilde v_0(r)=O(r^{-1})$ and
$\widetilde v_1(r)=O(r^{-2})$.
 For large $R$, the argument in the proof of \cite[Claim 5.16]{LawrieRodriguez19}, applied with \eqref{eq:projection-split}, \eqref{eq:dyadic-increments},
and \eqref{eq:second-channel}, yields
\[
 \left(\frac34\right)^j\bigl(|\widetilde v_0(R)|+|\widetilde v_1(R)|\bigr)
 \le |\widetilde v_0(2^jR)|+|\widetilde v_1(2^jR)|\lesssim  \left(\frac12\right)^j\frac{1}{R }.
\]
Thus both coefficients vanish for large $R$. Hence, \eqref{eq:projection-split} and \eqref{eq:second-channel} imply that
$\boldsymbol{\widetilde u}_*(0)\equiv 0$ on $r>R$.

\medskip\noindent
\underline{Step 4}: Rigidity.
We prove $\boldsymbol u_*(t)=\boldsymbol u_Q$ for all $t\in\bbR$.
The potential bound \eqref{eq:exterior-potential} allows us to use
\cite[Proposition~4.7]{DKM20Exterior} on the exterior region $\{r>R+|t|\}$, for each $R>0$.
Any nonzero compactly supported $\boldsymbol{\widetilde{u}}_*(0)$ produces
positive exterior energy in one time direction, which contradicts
\eqref{eq:difference-nonradiative}. Hence
$\boldsymbol u_*(0)=\boldsymbol{u}_{\ell}\in\mathcal H$.
The cutoff argument in the proof of Theorem~\ref{thm:gwp}
shows that the inverse image of $\boldsymbol{u}_{\ell}$ is a
stationary solution to \eqref{eq:energy-weak-form} across $r=0$. Therefore,
Theorem~\ref{thm:stationary} implies $\boldsymbol{u}_{\ell}=\boldsymbol{u}_Q$, which completes the proof of Proposition~\ref{prop:rigidity}.
\end{proof}
\subsection{Proof of soliton resolution}
\label{sec:transformed-completion}

We first prove Proposition~\ref{prop:transformed-resolution}. The return
to the original variable will then prove Theorem~\ref{thm:main}.

\begin{proof}[Proof of Proposition~\ref{prop:transformed-resolution}]
We show that every accumulation point of the remainder as $t\to\infty$
is the initial datum of a global solution to which Proposition~\ref{prop:rigidity}
applies. By Proposition~\ref{prop:compactness}, 
\[
 K=\overline{\{\boldsymbol{u}(t)-\boldsymbol{u}_L(t):t\ge0\}}^{\mathcal H}
\]
is compact. Let $\boldsymbol w=\boldsymbol u-\boldsymbol u_L$ and
choose a sequence $t_j\to\infty$. We may choose a subsequence,
still denoted by $t_j$, such that
\[
 \boldsymbol w(t_j)\to\boldsymbol u_{*,0}
 \quad\text{in }\mathcal H.
\]
Let $\boldsymbol u_*$ be the global solution in
Proposition~\ref{prop:transformed-gwp}, with
$\boldsymbol u_*(0)=\boldsymbol u_{*,0}$. For every $B,T>0$, the strong Huygens principle, continuous
dependence, and finite speed of propagation imply
\[
 \boldsymbol w(t_j+\cdot)\to\boldsymbol u_*
 \quad\text{in }C([-T,T];\mathcal H_B).
\] For each fixed $s$, compactness of $K$ and convergence in
$\mathcal H_B$ for every $B>0$ imply
\[
 \boldsymbol w(t_j+s)\to\boldsymbol u_*(s)
 \quad\text{in }\mathcal H.
\]
Hence $\{\boldsymbol u_*(s):s\in\bbR\}\subset K$, and so
Proposition~\ref{prop:rigidity} implies
$\boldsymbol u_{*,0}=\boldsymbol u_Q$. Thus every sequence $t_j\to\infty$ has a subsequence along which
$\boldsymbol u(t_j)-\boldsymbol u_L(t_j)\to\boldsymbol u_Q$
in $\mathcal H$, which proves \eqref{eq:transformed-resolution}.
\end{proof}

\begin{proof}[Proof of Theorem~\ref{thm:main}]
Let $\boldsymbol u_L$ be the free wave from
Proposition~\ref{prop:transformed-resolution}, and let
$\boldsymbol\psi_L=r\boldsymbol u_L$.
Then $\boldsymbol\psi_L$ is a finite-energy solution of
\eqref{eq:linear-angular}. For every $B>0$, we first show
\begin{equation}\label{eq:interior convergence}
 \boldsymbol\psi(t)-\boldsymbol Q\to0
 \quad\text{in }\mathcal E_{0,0}\text{ on }r<B.
\end{equation}
Let $\boldsymbol\psi_B(t)$ be the inverse image of
$\boldsymbol u_Q+\chi_B(\boldsymbol u(t)-\boldsymbol u_Q)$.
By \eqref{eq:transformed-resolution} and the strong Huygens principle,
\[
 \chi_B(\boldsymbol u(t)-\boldsymbol u_Q)
 \to0\quad\text{in }\mathcal H.
\]
Since $Q/r,Q'\in L^\infty$, subtraction of $\boldsymbol Q$
defines a continuous map from $\mathcal E_{0,n}$ to
$\mathcal E_{0,0}$. Proposition~\ref{prop:nonlinear-transform}
and the equality $\boldsymbol\psi_B(t)=\boldsymbol\psi(t)$ on $r<B$
therefore prove \eqref{eq:interior convergence}. Next, the transformation \eqref{eq:transform} yields
\[
 \psi- Q-\psi_L
 =r(u-u_Q-u_L)
 +\int_Q^\psi\left(1-\frac{h(r,\theta)}r\right)d\theta.
\]
We differentiate this identity and apply Hardy's inequality,
using $|h-r|\le r^{-1}$ and $|\partial_r(h/r)|\le2r^{-3}$.
Then \eqref{eq:uniform-energy} and
\eqref{eq:transformed-resolution} yield
\begin{equation*}
 \limsup_{t\to\infty}\norm{\boldsymbol\psi(t)-\boldsymbol Q-\boldsymbol\psi_L(t)}_{(\dot H^1\times L^2)(\bbR^3;\,r>R)}\lesssim_E R^{-2},\qquad R>1.
\end{equation*}
Together with \eqref{eq:interior convergence} and the strong
Huygens principle, this bound proves
\eqref{eq:main decomposition} as $R\to\infty$. 

Finally, \eqref{eq:interior convergence} and conservation of
energy imply \eqref{eq:excess-energy} since
\[
 \lim_{t\to+\infty}E(\boldsymbol\psi(t);r>R)
 =E(\boldsymbol\psi(0))-E(\boldsymbol Q;r<R)
 \qquad\text{for }R>0.
\]
This completes the proof of Theorem~\ref{thm:main}.
\end{proof}

\appendix
\section{Persistence of Sobolev regularity}\label{app:persistence}

We prove local persistence of Sobolev regularity using
Lemma~\ref{lem:lwp}, without assuming energy conservation.
This supplies the smooth approximations needed in
Proposition~\ref{prop:transformed-gwp}. Global persistence of regularity then follows
by applying the local statement on every compact interval. Let $\phi$
be the reference angle in
\eqref{eq:reference-angle}, and set $\boldsymbol{g}=(g,\dot g)$, where $\boldsymbol{g} = r^{-1} (\boldsymbol{\psi}-\boldsymbol{\phi})$.

\begin{prop}[Persistence of Sobolev regularity]\label{prop:persistence}
(i) For $s>5/2$, the map
$\boldsymbol g\mapsto\boldsymbol u$ induced by
\eqref{eq:transform} is a homeomorphism of
$H^s(\bbR^5)\times H^{s-1}(\bbR^5)$.

(ii) Let $s\ge1$ and
$\boldsymbol u_0\in H^s(\bbR^5)\times H^{s-1}(\bbR^5)$.
Let $\boldsymbol u\in C(I;\mathcal H)$ be the solution in
Lemma~\ref{lem:lwp}, with $\boldsymbol u(0)=\boldsymbol u_0$.
Then
\begin{equation}
\boldsymbol u\in C\bigl(I;H^s(\bbR^5)\times H^{s-1}(\bbR^5)\bigr).\label{eq:gwp-hs}
\end{equation}
The corresponding continuous dependence holds on compact time intervals.
\end{prop}

We estimate derivatives of $F(r,u)$ for $0\le r\le1$, where $\phi=0$. We have
\[
 u=A(r,g),\qquad
 A(r,g)=\int_0^g\sqrt{1+2\left({\sin(r\eta)}/r\right)^2}d\eta.
\]
In $F_r$, $F_{rr}$, and $F_{ru}$, derivatives with respect to $r$ are taken with $u$ fixed.

\begin{lem}[Derivatives of $F$]\label{lem:gwp-coefficients}
Let $0<\delta<\frac{\pi}{2}$. If $0<r\le1$ and
$|H^{-1}(r,r^2u)|\le\delta$, then
\begin{equation}\label{eq:gwp-explicit}
|F_u|\lesssim_\delta |u|,\; |F_{uu}|\lesssim_\delta 1,\;
 |F_r|\lesssim_\delta r|u|^3,\; |F_{ru}|\lesssim_\delta r|u|^2,\;
 |F_{rr}|\lesssim_\delta |u|^3.
\end{equation}
\end{lem}
\begin{proof}
All implicit constants depend only on $\delta$. By oddness, we may assume $u\ge0$. The first bound in \eqref{eq:gwp-explicit} follows from
Lemma~\ref{lem:nonlinear-estimates}.
Since $b=\sin^2(rg)/r^2$, implicit differentiation of $u=A(r,g)$ yields
\[
 1+b\sim 1+u,\qquad |b_u|\le\sqrt2,\qquad
 |b_r|\lesssim ru^2,\qquad |b_{rr}|\lesssim u^2.
\]
The identity $F_u=3b+b/(1+2b)^2$ and the chain rule then imply
\[
 |F_{uu}|\lesssim1,\qquad |F_{ru}|\lesssim ru^2,\qquad
 |F_{rru}|\lesssim u^2.
\]
We obtain \eqref{eq:gwp-explicit} by integration and monotonicity.
\end{proof}

\begin{proof}[Sketch of proof of Proposition~\ref{prop:persistence}]
$ $

\noindent
\underline{Step 1}. We establish the regularity correspondence and reduce continuation
to an interior estimate.
For $s>5/2$, the Sobolev embedding
$H^s(\bbR^5)\hookrightarrow L^\infty(\bbR^5)$
and the composition and product inequalities, applied to
$u=A(r,g)$ and its inverse, identify the
$H^s\times H^{s-1}$ regularity of $\boldsymbol g$ and
$\boldsymbol u$ near $r=0$.
Away from $r=0$, multiplication by $r$ identifies
$H^s\times H^{s-1}$ regularity in $\bbR^5$ and $\bbR^3$.
Sobolev embedding and the same composition and product estimates
yield continuity in both directions.

Let $[0,T]\subset I$. Continuity in $\mathcal H$, integration of $u_t$, and
Lemma~\ref{lem:lwp} imply
\[
 \sup_{0\le t\le T}\|\boldsymbol u(t)\|_{H^1\times L^2}+\|u\|_{L^4_{[0,T]}L^4_x}\le M.
\]
The radial Sobolev inequality implies
$|r^2u(t,r)|\lesssim M\sqrt r$.
For sufficiently small $a>0$ with $4a<1$, monotonicity of
$H(r,\cdot)$ therefore implies
$|\psi(t,r)|\le\delta<\pi/2$ on $r\le4a$. On $r\ge a/4$, equation~\eqref{eq:relative-equation} becomes
\begin{equation}\label{eq:gwp-exterior}
 \partial_t^2(ru)-\partial_r^2(ru)-\frac2r\partial_r(ru)=rF(r,u)-\frac{2u}{r}.
\end{equation}
The right-hand side of \eqref{eq:gwp-exterior} is smooth in $r$
and $ru$. Since $|ru|\lesssim_a M$ by the radial Sobolev inequality,
its mixed derivatives are uniformly bounded for $r\ge a/4$. To analyze the solution near the origin, we localize:
\[
 \partial_t^2(\chi_a u)-\Delta_5(\chi_a u)=\chi_{2a}F(r,\chi_a u)+G,
\]
where $G$ is supported on $a\le r\le2a$. Let $\sigma\ge 3$ and $J=[t_0,t_1]\subset[0,T]$ with $|J|\le a/4$.
For $\boldsymbol u\in C(J;H^\sigma\times H^{\sigma-1})$,
finite speed of propagation, the Duhamel formula, and
Moser's and Gronwall's inequalities for
\eqref{eq:gwp-exterior} imply
\begin{equation}\label{eq:gwp-exterior-regularity}
 \sup_{t\in J}\|\boldsymbol u(t)\|_{(H^\sigma\times H^{\sigma-1})(r>a/2)}\!+\!\|G\|_{L^\infty_JH^{\sigma-1}}\lesssim_{M,a,\sigma}1+\|\boldsymbol u(t_0)\|_{H^\sigma\times H^{\sigma-1}}.
\end{equation}
\underline{Step 2}. For $H^3\times H^2$ initial data and sufficiently short $J$, we prove
\begin{equation}\label{eq:gwp-h3-tame}
 \|\boldsymbol u\|_{C(J;H^3\times H^2)}\lesssim_{M,a}1+\|\boldsymbol u(t_0)\|_{H^3\times H^2}.
\end{equation}
The Gagliardo--Nirenberg and the radial Sobolev inequalities imply
\[
 \|\nabla(\chi_a u)\|_4^2\lesssim\|\chi_a u\|_4\|\nabla^2(\chi_a u)\|_4,\quad r|\chi_a u|^{2/3}\lesssim M^{2/3}.
\]
H\"older's inequality, Sobolev embeddings, and interpolation yield
\[
\begin{aligned}
 \|r|\chi_a u|^2\nabla(\chi_a u)\|_2
 &\lesssim M^{2/3}\|\chi_a u\|_{10/3}^{4/3}
             \|\nabla(\chi_a u)\|_{10}
 \lesssim M^2\|\chi_a u\|_{H^3},\\
 \|\chi_a u\|_6^3
 &\lesssim\|\chi_a u\|_{\dot H^1}^2
             \|\chi_a u\|_{\dot H^3}
 \lesssim M^2\|\chi_a u\|_{H^3}.
\end{aligned}
\]
Lemma~\ref{lem:gwp-coefficients} consequently implies
\[
 \|\chi_{2a}F(r,\chi_a u)\|_{L^1_JH^2}\lesssim |J|^{1/2}\|\chi_a u\|_{L^4_JL^4}\|\nabla^2(\chi_a u)\|_{L^4_JL^4}+C_M|J|\|\chi_a u\|_{L^\infty_JH^3}.
\]
On $2a<r<4a$, we have $\chi_a u=0$ and $F(r,0)=0$.
Hence all terms containing derivatives of $\chi_{2a}$ vanish.
Lemma~\ref{lem:lwp} bounds
$\|\chi_a u\|_{L^4_JL^4}$ uniformly on $[0,T]$.
For $|J|$ sufficiently small depending only on $M,a$, the Strichartz inequality
\eqref{eq:linear}, applied to $\nabla^2(\chi_a u)$,
together with \eqref{eq:gwp-exterior-regularity},
implies \eqref{eq:gwp-h3-tame}.
Iteration on finitely many intervals yields
$\boldsymbol u\in C([0,T];H^3\times H^2)$.

\noindent
\underline{Step 3}. We now prove \eqref{eq:gwp-hs}. Lemma~\ref{lem:lwp} and integration of $u_t$ imply local
Lipschitz continuity of the solution map in $H^1\times L^2$.
Nonlinear interpolation of the solution map with
\eqref{eq:gwp-h3-tame} proves persistence and continuous
dependence for $1<s<3$. For $s\ge 3$, \eqref{eq:gwp-h3-tame} and the Sobolev embedding
$H^3(\bbR^5)\hookrightarrow L^\infty(\bbR^5)$
control $\|\chi_a u(t)\|_\infty$ on $[0,T]$.
Moser's inequality implies
\[
 \|\chi_{2a}F(r,\chi_a u)\|_{H^{s-1}}\lesssim1+\|\chi_a u\|_{H^{s-1}}.
\]
The Duhamel formula, \eqref{eq:gwp-exterior-regularity} with
$\sigma=s$, and Gronwall's inequality prove \eqref{eq:gwp-hs}.
Continuous dependence follows from the local
$H^s\times H^{s-1}$ theory.
\end{proof}
\bibliographystyle{abbrv}
\bibliography{Skyrme}

\end{document}